\documentclass[sn-mathphys-ay]{sn-jnl}% Math and Physical Sciences Author Year Reference Style
\usepackage{graphicx}%
\usepackage{multirow}%
\usepackage{amsmath,amssymb,amsfonts}%
\usepackage{amsthm}%
\usepackage{mathrsfs}%
\usepackage[title]{appendix}%
\usepackage{xcolor}%
\usepackage{textcomp}%
\usepackage{manyfoot}%
\usepackage{booktabs}%
\usepackage{algorithm}%
\usepackage{algorithmicx}%
\usepackage{algpseudocode}%
\usepackage{listings}%
\begin{document}

\title[Article Title]{Modeling and solving an integrated periodic vehicle routing and capacitated facility location problem in the context of solid waste collection}


\author*[1]{\fnm{Gonz\'alez} \sur{Bego\~na}}\email{bego.landin@ulpgc.es}

\author*[2]{\fnm{Rossit} \sur{Diego}}\email{diego.rossit@uns.edu.ar}

\author[3]{\fnm{Frutos} \sur{Mariano}}\email{mfrutos@uns.edu.ar}

\author[1]{\fnm{M\'endez} \sur{M\'aximo}}\email{maximo.mendez@ulpgc.es}

\affil*[1]{\orgdiv{Instituto Universitario SIANI}, \orgname{Universidad de Las Palmas de Gran Canaria}, \orgaddress{\street{Parque Cient\'ifico y Tecnol\'ogico, Campus Universitario de Tafira}, \city{Las Palmas de Gran Canaria}, \postcode{35017}, \state{Canarias}, \country{Spain}}}

\affil[2]{\orgdiv{INMABB, Departament of Engineering}, \orgname{Universidad Nacional del Sur-CONICET}, \orgaddress{\street{Alem Av. 1253}, \city{Bah\'ia Blanca}, \postcode{8000}, \state{}, \country{Argentina}}}

\affil[3]{\orgdiv{IIESS, Departament of Engineering}, \orgname{Universidad Nacional del Sur-CONICET}, \orgaddress{\street{Alem Av. 1253}, \city{Bah\'ia Blanca}, \postcode{8000}, \state{}, \country{Argentina}}}


\abstract{Few activities are as crucial in urban environments as waste management. Mismanagement of waste can cause significant economic, social, and environmental damage. However, waste management is often a complex system to manage and therefore where computational decision-support tools can play a pivotal role in assisting managers to make faster and better decisions. In this sense, this article proposes, on the one hand, a unified optimization model to address two common waste management system optimization problem: the determination of the capacity of waste bins in the collection network and the design and scheduling of collection routes. The integration of these two problems is not usual in the literature since each of them separately is already a major computational challenge. On the other hand, two improved exact formulations based on mathematical programming and a genetic algorithm (GA) are provided to solve this proposed unified optimization model. It should be noted that the GA considers a mixed chromosome representation of the solutions combining binary and integer alleles, in order to solve realistic instances of this complex problem. Also, different genetic operators have been tested to study which combination of them obtained better results in execution times on the order of that of the exact solvers. The obtained results show that the proposed GA is able to match the results of exact solvers on small instances and, in addition, can obtain feasible solutions on large instances, where exact formulations are not applicable, in reasonable computation times.}




\maketitle

\section{Introduction}\label{sec:introduction}

%The Introduction section, of referenced text \cite{bib1} expands on the background of the work (some overlap with the Abstract is acceptable). The introduction should not include subheadings.
%
%Springer Nature does not impose a strict layout as standard however authors are advised to check the individual requirements for the journal they are planning to submit to as there may be journal-level preferences. When preparing your text please also be aware that some stylistic choices are not supported in full text XML (publication version), including coloured font. These will not be replicated in the typeset article if it is accepted. 

Today, cities face numerous challenges to ensure the livability of the urban environment and the well-being of their citizens. One critical aspect that demands accurate handling is the proper collection of waste generated as part of the various city activities. A primary source of this waste is households, which contribute to what is known as Municipal Solid Waste (MSW). MSW is known for its heterogeneity, as it comprises different types of materials and is characterized by a distributed generation~\citep{rossit2017application}. In contrast to other waste types, such as industrial or sanitary waste, which are typically produced in large quantities at specific locations, MSW is generated at numerous city locations in smaller amounts. This poses a logistical challenge for authorities tasked with collecting this waste, often on a daily or two-day basis, from several scattered sources.

When it comes to collecting waste from households, various systems can be implemented. One traditional approach is the door-to-door system, where waste is collected individually from each household. This method, still prevalent in many cities, involves the collection vehicle visiting every household. Another common method is the community bin system, in which citizens deposit waste in bins distributed throughout the city. The points where these community bins are located, commonly referred to as \textit{collection points}, are visited, not necessarily every day, by collection vehicles that empty the bins. Recent research suggests that the use of a point-based system is preferable to door-to-door collection in crowded places, as it can enhance the cost-efficiency of the logistics of the system, for example, by minimizing the distances that waste collection trucks have to travel~\citep{rossi2022comparison,rossit2022waste}. This logistic cost reduction may be especially notable in nations with already elevated logistics costs, such as Argentina~\citep{alalog2023}, and in particular for the city of Bah\'ia Blanca, which exceeds the national average in waste management expenditure~\citep{rossit2023systemwaste}, and from which the case study presented in this paper were taken.

The establishment of the collection point network in the community bin system involves determining the capacity of these locations. Thus, the decision regarding the capacity of collection points is closely interconnected with the facet of the MSW system that comes immediately after in the reverse logistics chain, i.e., the waste collection. The management of MSW collection includes the planning of routes, i.e., the actual paths travelled within city streets to collect the waste deposited at these collection points, and the arrangement of these routes within a planning timeframe to create a regular, usually weekly, schedule. These two problems: the determination of collection point capacity and the design of collection routes, are interrelated. Thus, given a waste generation rate, the storage capacity of the collection points will dictate the frequency of collection truck visits to prevent overflow~\citep{maheo2020benders}. However, both problems are usually solved separately due to the computational complexity that an integrated approach would imply. In fact, only solving the point-based network design involves solving a capacitated facility location problem, which is known to be an NP-hard problem~\citep{cornuejols1991comparison}. In developing countries, these issues are often addressed based on practitioner experience rather than computer-aided tools~\citep{maalouf2023waste}, including Argentina~\citep{cavallin2020application}. While practical experience is valuable, the academic literature offers insights into cost minimization and environmental preservation through computer-aided decision support methodologies.

Thus, the aim of this study is to propose an integrated approach to simultaneously address the design of both network and collection routes for MSW, which is an area were there is still a vacancy in the related literature. The remainder of the paper is organized as follows: Section~\ref{sec:lit_review} provides a review of the related literature and outlines the main contributions of this work to the state-of-the-art in the subject. Section~\ref{sec:problem_description} presents both the conceptual description and the mathematical formulation of the problem addressed. Section~\ref{sec:resolution} describes the resolution approaches considered to solve the proposed mathematical model. Section~\ref{sec:GAparam} reports the instances used in the paper as well as the $4\times3^3$ factorial design conducted for tuning the GA. Section~\ref{sec:NumResults} shows the numerical results achieved with both resolution approaches and performs a comparative analysis of them. Finally, Section~\ref{sec:conclusions} presents the conclusions and future work.

\section{Literature review}\label{sec:lit_review}

As mentioned above, the integrated problem of simultaneously determining the capacity of collection points and designing collection routes has received limited attention in the literature. \cite{hemmelmayr2014models} proposed the first model based on mathematical programming that simultaneously determines the capacity of waste bins alongside waste collection routes and their scheduling. They proposed different conceptual models of the bin location part considering whether there is an initial configuration of bins to be relocated or not, considering whether there are space limits, and considering whether there is waste source classification. They solved small synthetic instances with an exact solver and complemented it with a Large Neighborhood Search (LNS) algorithm for solving real-world instances of an unspecified city in Italy. 

In the context of recyclable waste, \cite{hemmelmayr2017periodic} presented a location-routing approach that aims to simultaneously optimize the location of collection points, the capacity of bins to be installed within each collection point and the routing schedule. They proposed a mathematical formulation for solving small instances and an adaptive LNS for larger instances, which is able to obtain near-optimal solutions in small instances. The largest instance they addressed considered a maximum of 50 collection points. Also in the context of recyclable waste, \cite{Cubillos2020} proposed a similar model that was addressed with a hybrid Variable Neighborhood Search (VNS) algorithm. In this case, it should be noted that the waste collection problem was limited to a single day, so it did not involve periodic routing, and that the authors chose an uncapacitated vehicle for waste collection, so it could collect all the garbage in a single trip. All of this is a major simplification that reduces the complexity of the problem. 

Subsequently, \cite{glaser2021introduction} proposed a model that extended the problem of~\cite{hemmelmayr2014models} by introducing predefined alternatives for bin collection frequencies and considering user allocation to collection points. They presented an adaptive LNS heuristic to solve the problem, which could outperform the results of plain CPLEX.

In a related context, \cite{roy2022iot} presented a model that tackles waste bin allocation and vehicle routing using VNS with Ant Colony Optimization, incorporating smart bins with alert capability for imminent filling. They analyze the presented model with realistic instances of Seoul (Korea) with 15 collection points, concluding that the model is effective in finding good quality solutions in a reasonable time, however, it should be noted that their model does not incorporate decision making on what type and capacity of bins to place at each collection point. 

Lastly, \cite{maheo2023solving} presented integrated models to solve a similar allocation-routing problem, employing a novel resolution method based on Benders decomposition and valid cuts. The proposed method is able to outperform plain CPLEX in realistic instances, with at most 7 collection points. As in previous works, the authors considered predefined alternatives for bin collection frequencies.

Recently, other studies have addressed similar location-routing problems within the reverse logistic chain of municipal waste. For instance, \cite{Han2024Optimizing} focused on a location-routing problem in which waste generation is grouped by residential sectors. Their model integrates decisions such as transporting waste from residential sectors to transfer stations using small vehicles, and then collecting waste at transfer stations and transporting it to disposal centers using large vehicles. However, their proposed model, which aims to minimize the distance traveled by both types of vehicles, does not take into account either the capacity restrictions of the facilities or a periodic schedule for routing. To solve this problem, the authors presented a mathematical formulation and used an approach based on the Dantzig-Wolfe decomposition, which demonstrated superior performance compared to plain CPLEX in real-world instances of Nanchang, China. For a similar problem but considering a stochastic waste generation rate following a known normal distribution, \cite{niu2024multi} presented a multi-objective model to simultaneously optimize cost, greenhouse gas emission and citizens' satisfaction. However, neither the capacity of the facilities nor a periodic schedule of the routing were considered. To solve the optimization problem, they proposed an evolutionary algorithm hybridized with a decision tree, which is able to outperform other classical multi-objective evolutionary algorithms in real-world instances of the Chinese city of Beijing.

\subsection{Contribution of this work}

After a thorough review of the literature related to MSW collection by the community bin system, it has been observed that the integrated problem of simultaneously determining the capacity of collection points and designing collection routes has received limited attention. In this framework, the present document contributes to the state-of-the-art with the following items:

\begin{enumerate}
	\item \textit{Expansion of the conceptual model}. The conceptual model considered in other works is extended by including decisions on the capacity of the facilities and avoiding predefined alternatives for MSW collection frequencies. This approach, although it implies an increase in computational complexity, allows greater flexibility both in choosing the bin combination to be placed at each collection point and in establishing collection frequency of MSW accumulated in each one of them.
	\item \textit{Novel mathematical formulation}. A new mathematical formulation is proposed to formalize the extended problem that began to be developed in~\cite{Rossit2024}, although the model proposed here is substantially different since new features have been added. Now the result is a Mixed Integer Quadratic Programming (MIQP) model and a different linearized Mixed Integer Linear Programming (MILP) model, which are used with exact solvers to address the problem.
	\item \textit{Proposal of a mixed GA}. Since it is widely recognized that metaheuristic approaches are suitable for solving computationally complex optimization problems, a GA with a mixed binary-permutation encoding is proposed to solve the target problem. 
	\item \textit{Computational experimentation with real-world instances}. We start by comparing the performance of the two resolution approaches considered (exact solver versus AG) using relatively small instances, in line with other previous studies using exact approaches, since this type of approach is usually only applicable for not very large instances. We end with an evaluation of the ability of the AG to solve real-world instances using one with 163 collection points corresponding to an area of the Argentinean city of Bah\'ia Blanca. It should be noted that the size of this instance is considerable compared to the instances considered in other works reviewed in the related literature.
\end{enumerate}

\section{The integrated periodic vehicle routing and capacitated facility location problem in the context of MSW collection}\label{sec:problem_description}
The optimization problem addressed in this section arises by simultaneously considering two different aspects of the MSW reverse logistics chain, namely:

\begin{enumerate}
	\item the waste bin allocation problem to set the types and number of bins to be installed in predefined locations, according to the local waste generation, and the visit schedule of the collection vehicle. 
	\item the waste collection problem to define the collection points to be visited each day, the visit order of these points and the scheduling of the daily routes in the planning horizon (since not all collection points must to be visited every day). This corresponds to a Periodic Capacitated Vehicle Routing Problem (PCVRP)~\citep{beltrami74}.
\end{enumerate}

This joint approach represents an inventory-routing problem~\citep{archetti2022comparison} with the additional feature that the storage capacity of each collection point is not predefined, instead it must be determined within the optimization process to minimize both the cost of the routing schedule and the cost of bin installation and maintenance. Additionally, realistic constraints are taken into account both for bin allocation and for the design of the daily routes in the planning horizon considered. These constraints include: a discrete number of capacity alternatives for bin combinations at collection points derived from the type of commercial bins available, a limit on the waste accumulation at collection points to avoid exceeding the capacity of the bins installed, a restriction on the amount of waste collected by each collection vehicle on each route, to ensure that it does not exceed its capacity, a limitation on the number of vehicles required per day that not exceed the fleet size, and the guarantee that the time required to complete each route from the departure of the collection vehicle from the depot to its return, including service times at each collection point and unloading time at the depot, does not exceed the drivers' work shift.

The mathematical model proposed extends the work of \cite{maheo2020benders,maheo2023solving} in that it does not simply consider a predefined set of feasible visit combinations to collect waste from collection points, but rather any possible visit combination can be used as long as the collection point sites do not overflow (i.e., as long as the waste accumulated at each collection point does not exceed the capacity of the bin combination chosen for that point). This is a feature that differentiates it from other models proposed in relevant papers addressing similar problems in the literature~\citep{glaser2021introduction,hemmelmayr2014models}. In addition, the time of service at each collection point is not fixed as in ~\cite{maheo2020benders,maheo2023solving}, but depends on the type and number of bins installed.

Let $n_I$, $n_B$, $n_V$ and $n_T$ be the number of collection points, bin combinations, collection vehicles and days in the planning horizon, respectively. A collection point, $i\in I = \{1, 2, ..., n_I\}$, is a predefined location in an urban area where waste bins can be installed. We define $I^0 =  I\cup \{0\}$, where $0$ represents the depot where the collected waste is deposited and where each vehicle, $v \in V = \{1, 2, ..., n_V\}$, starts and ends its daily routes. A bin combination, $b \in B = \{1, 2, ..., n_B\}$, is a set of bins that can be installed at a collection point, taking into account physical space constraints. Therefore, each element of $B$ can consist of a single bin or a combination of bins. The method used to determine the bin combinations considered in this paper can be found in~\cite{maheo2023solving}. Lastly, subset $T' \subset T = \{1, 2, ..., n_T\}$ is defined as the days of the planning horizon on which waste collection is not carried out due to drivers' rest days, for example, Sundays in the city of Bah\'ia Blanca, Argentina

The following parameters are also defined:

\begin{itemize}
	\item[] $Q$: Collection vehicle capacity ($m^3$).
	\item[] $C_{ij}$: Travel time (minutes) between points $i\in I^0$ and $j\in I^0$.
	\item[] $S_b$: Service time (minutes) of bin combination $b\in B$.
	\item[] $T_U$: Unloading time (minutes) of collection vehicles at the depot.
	\item[] $W_i$: Amount of waste deposited daily ($m^3/day$) at collection point $i\in I$.
	\item[] $CAP_b$: Capacity ($m^3$) of bin combination $b\in B$.
	\item[] $CIN_b$: Adjusted cost over the planning horizon (US\$) of bin combination $b\in B$, taking into account both acquisition and installation costs, as both are to be amortized over the useful life of the bins.
	\item[] $C_{CV}$: Cost of collection vehicles (US\$/min).
	\item[] $T_{L}$: Length of working day for drivers (minutes).
\end{itemize}

%It should be noted that $CIN_b$ is an adjusted cost. This is because the model considers two different decision and cost levels: i) a strategic decision-making level involving the purchase and installation of bins that are expected to last for several years; and ii) a tactical decision-making level involving the transportation costs in the routing plan and the maintenance of the bins. Therefore, the cost of purchasing and installing the bins needs to be amortized over their useful life.

Waste collection is assumed to take place at the end of the day, once all waste generated (daily) has been deposited at the corresponding collection points. With this in mind, the following variables are also defined:

\begin{itemize}
	\item[] $x_{ijvt}$: 1 if vehicle $v\in V$ travels from point $i\in I^0$ to point $j\in I^0$ on day $t\in T$, and 0 otherwise.  
	\item[] $y_{ijvt}$: load of vehicle $v\in V$ when it travels from point $i\in I^0$ to point $j\in I^0$ on day $t\in T$.
	\item[] $w_{it}$: amount of waste accumulated at collection point $i\in I$ at the end of day $t\in T$.
	\item[] $w_i^{max}$: maximum daily amount of waste accumulated at collection point $i\in I$ in the planning horizon.
	\item[] $n_{bi}$: 1 if bin combination $b\in B$ is selected for collection point $i\in I$, and 0 otherwise.
\end{itemize}

Finally, Eq.~\eqref{eq:eq_TT} calculates $TT_{vt}$, the total service time of vehicle $v \in V$, on day $t \in T$, from the time it leaves the depot until it returns to it, including the service times at each collection point visited and the unloading time at the depot.

\begin{equation} 
	\label{eq:eq_TT}
	TT_{vt} = T_U \sum_{i\in {I}} x_{i0vt} + \sum_{i\in {I^0}} \sum_{j\in{I^0}} x_{ijvt} \left( C_{ij}  + \sum_{b\in B} S_{b} n_{bj} \right)
\end{equation}

Taking all these elements into account, the following mathematical model is proposed:

\begin{subequations} 
	\label{eq:eq_2}
	Minimize
	\begin{equation} 
		\sum_{b\in B} \left( CIN_{b} \sum_{i\in I} n_{bi} \right) + C_{CV} \sum_{t\in T} \sum_{v\in V} TT_{vt}
		\tag{\ref{eq:eq_2}}
	\end{equation}
	\quad subject to
	\allowdisplaybreaks
	\begin{align}
		n_{b0} = 0,
		{\begin{array}{r}\forall b\in B \end{array}}\label{eq:eq_2a}\\[0.05in]
		\sum_{b\in B} n_{bi} = 1, 
		{\begin{array}{r}\forall i\in I \end{array}}\label{eq:eq_2b}\\[0.05in]
		\sum_{b\in B} CAP_{b} n_{bi} \geq w_i^{max}, 
		{\begin{array}{l}\forall i\in I \end{array}}\label{eq:eq_2c} \\[0.05in]
		x_{iivt = 0 \, 
			{\begin{array}{l}\forall %(i,v,t)\in I^0\times V\times T 
					i\in I^0, v\in V, t\in T
		\end{array}}}\label{eq:eq_2d}\\[0.05in]
		x_{ijvt} = 0 \, 
		{\begin{array}{l}\forall 
				%(i,j,v,t)\in I^0 \times I^0\times V\times T'
				i,j \in I^0, v\in V, t\in T'
		\end{array}}\label{eq:eq_2e}\\[0.05in]	
		\sum_{i\in I^0} x_{ijvt}- \sum_{i\in I^0} x_{jivt} = 0,\, 
		{\begin{array}{l}
				%\forall (j,v,t)\in I^0\times V\times T
				j \in I^0, v \in V, t \in T
		\end{array}} \label{eq:eq_2f}\\[0.05in]
		\sum_{i\in I} x_{0ivt} \leq 1,\,
		{\begin{array}{l}\forall v\in V, t\in (T-T') \end{array}}\label{eq:eq_2g}\\[0.05in]        
		TT_{vt} \leq T_{L},\,
		{\begin{array}{l}\forall 
				v\in V, t\in (T-T')
		\end{array}}       
		\label{eq:eq_2h}\\[0.05in]
		y_{ijvt} \leq Q x_{ijvt},\,
		{\begin{array}{l}\forall 
				i,j\in I^0, v\in V, t\in T
		\end{array}}\label{eq:eq_2i}\\[0.05in]
		\sum_{i\in I^0} y_{ijvt} + w_{jt} \leq \sum_{i\in I^0} y_{jivt} + Q \left( 1 - \sum_{i\in I^0} x_{ijvt} \right),\,\nonumber\\[0.05in]
		{\begin{array}{l}\forall 
				j \in I, v\in V, t\in T \end{array}}\label{eq:eq_2j}\\[0.05in]
		w_{it} = W_i + w_{i(t-1)} \left( 1 - \sum_{j\in I^0} \sum_{v\in V} x_{ijv(t-1)} \right), \nonumber \\
		{\begin{array}{l}\forall 
				i \in I, t\in (T-\{1\})
		\end{array}}\label{eq:eq_2k}\\[0.05in]
		w_{i1} = W_i + w_{in_T} \left( 1 - \sum_{j\in I^0} \sum_{v\in V} x_{ijv n_T} \right),
		{\begin{array}{l}\forall i\in I \end{array}}\label{eq:eq_2l}\\[0.05in]
		w_{it} \leq w_{i}^{max}, {\begin{array}{l}\forall i \in I, t\in T \end{array}}\label{eq:eq_2m}
	\end{align}
\end{subequations}

where (in matrix notation) 

\begin{itemize}
	\item[] $x = [x_{ijvt}] \in \mathcal{M}_{(n_I+1) \times (n_I+1) \times n_V \times n_T}\left(\{0,1\}\right)$, 
	\item[] $y = [y_{ijvt}] \in \mathcal{M}_{(n_I+1) \times (n_I+1) \times n_V \times n_T}\left(\mathbb{R}^{+} \cup \{0\}\right)$,    
	\item[] $w = [w_{it}] \in \mathcal{M}_{n_I \times n_T}\left(\mathbb{R}^{+} \cup \{0\}\right)$, 
	\item[] $w^{max} = [w_i^{max}] \in \mathcal{M}_{n_I}\left(\mathbb{R}^{+} \cup \{0\}\right)$,
	\item[] $n = [n_{bi}] \in \mathcal{M}_{n_B \times (n_I+1)}\left(\{0,1\}\right)$.
	\item[]
\end{itemize}

The mathematical model provided simultaneously addresses the problems of MSW collection and bin combination selection. Thus, the objective function~\eqref{eq:eq_2} calculates the overall cost of installation and maintenance of bins and MSW collection. The latter based on the total service time of the vehicles~\eqref{eq:eq_TT}. Regarding the restrictions, Eq.~\eqref{eq:eq_2a} ensures that the model does not select a bin combination for the depot while Eq.~\eqref{eq:eq_2b} states that only one bin combination can be placed at each collection point. Eq.~\eqref{eq:eq_2c} establishes that the storage capacity installed at each collection point has to be sufficient to contain the maximum quantity of waste accumulated at that point over the planning horizon. Eq.~\eqref{eq:eq_2d} shows that vehicles always move between different collection points and, once they have left the depot, they must visit at least one collection point before returning to the depot. Eq.~\eqref{eq:eq_2e} shows that vehicles are out of service on drivers' rest days. Eq.~\eqref{eq:eq_2f} states that each vehicle has to leave each visited collection point. Eq.~\eqref{eq:eq_2g} establishes that each vehicle may perform one or no routes per day. Therefore, the number of daily routes cannot exceed the size of the vehicle fleet. Eq.~\eqref{eq:eq_2h} prevents that each route, from the time a collection vehicle leaves the depot until it returns to the depot, including service times at each collection point visited and unloading time at the depot, can be performed in the working day. Eq.~\eqref{eq:eq_2i} forces not to exceed vehicle capacity on collection routes. Eq.~\eqref{eq:eq_2j} is a subroute elimination restriction that also keeps track of vehicle loading by updating variable $y_{ijvt}$~\citep{toth2014vehicle}. Eqs.~\eqref{eq:eq_2k} and~\eqref{eq:eq_2l} keep track of the waste accumulated at each collection point at the end of each day, which depends on the waste generated on the day and whether any collection vehicle has visited that point on the previous day. For example, in Eq.~\eqref{eq:eq_2k}, if a collection point $i$ is not visited the previous day, $t-1$, then $\sum_{j\in I^0} \sum_{v\in V} x_{ijv(t-1)} = 0$ and Eq.~\eqref{eq:eq_2k} is reduced to $w_{it} = W_i + w_{i(t-1)}$. Conversely, when a collection point $i$ is visited the previous day, $t-1$, then $\sum_{j\in I^0} \sum_{v\in V} x_{ijv(t-1)} = 1$ and Eq.~\eqref{eq:eq_2k} is reduced to $w_{it} = W_i$. Specifically, Eq.~\eqref{eq:eq_2l} defines the cyclic characteristic of planning since, on the first day of the planning horizon, the residuals accumulated on the last day of the planning horizon has to be considered, in order for the solution to be feasible when the problem is implemented as a PCVRP. Eq.~\eqref{eq:eq_2m} sets, for each collection point, the maximum amount of waste accumulated each day of the planning horizon, in order to estimate the capacity of the bin combination to be installed at each collection point (linked to Eq.~\eqref{eq:eq_2c}).

\section{Resolution approaches}\label{sec:resolution}

Two different approaches have been chosen to solve the problem described in Section~\ref{sec:problem_description}. Specifically, two exact methods and a genetic algorithm have been considered.

\subsection{Exact methods}

The mathematical model described in Section~\ref{sec:problem_description} is a MIQP formulation due to the presence of a bilinear term both in the expression to compute the travel time $TT_{vt}$~\eqref{eq:eq_TT} which affects the objective function~\eqref{eq:eq_2} and the time limit restriction~\eqref{eq:eq_2h}, and in restrictions \eqref{eq:eq_2k} and \eqref{eq:eq_2l}, although the latter two can be replaced by~Eqs.~\eqref{eq:eq_l1}--\eqref{eq:eq_l3} which no longer have bilinear terms. 

\begin{subequations}
	\allowdisplaybreaks
	\begin{alignat}{1}
		w_{it} \geq W_i + w_{i(t-1)} - BigM \sum_{j\in I^0} \sum_{v\in V} x_{ijv(t-1)},  
		{\begin{array}{l}\forall i \in I, t\in (T-\{1\}) \end{array}}\label{eq:eq_l1}\\[0.05in]
		w_{i1} \geq W_i + w_{in_T} - BigM \sum_{j\in I^0} \sum_{v\in V} x_{ijvn_T},
		{\begin{array}{l}\forall i\in I \end{array}}\label{eq:eq_l2}\\[0.05in]
		w_{it} \geq W_i, {\begin{array}{l}\forall i\in I, \forall t\in T \end{array}}\label{eq:eq_l3}
	\end{alignat}
\end{subequations}

{\noindent}where $BigM$ is a sufficiently large value that resets the accumulated waste ($w_{it}$) at the collection point $i$ when visited. For example, in Eq.~\eqref{eq:eq_l1}, if the collection point $i$ is not visited the previous day, $t-1$, $\sum_{j\in I^0} \sum_{v\in V} x_{ijv(t-1)} = 0$ and Eq.~\eqref{eq:eq_l1} is reduced to $w_{it} \geq W_i + w_{i(t-1)}$. Conversely, when the collection point is visited the previous day, $t-1$, $\sum_{j\in I^0} \sum_{v\in V} x_{ijv(t-1)} = 1$ and Eq.~\eqref{eq:eq_l1} is reduced to $w_{it} \geq 0$, since it is a non-negative continuous variable. This is why it is necessary to add Eq.~\eqref{eq:eq_l3}. Note that, in either case, the value of $w_{it}$ will be pushed to the minimum ($w_{it} = W_i + w_{i(t-1)}$ or $w_{it} = W_i$, respectively) due to the sense of the optimization. A similar behaviour is observed in Eq.~\eqref{eq:eq_l2}. After careful consideration, $BigM$ was set as $BigM = \max\limits_{{b\in B}} \{CAP_b\}$. 

Usually, linear models have certain advantages over quadratic models, which tend to be more mature than their non-linear counterparts~\citep{rodriguez2013comparative}. Taking this into account, a MILP model was also implemented by applying the linearization technique proposed by \cite{glover1975improved}. Thus, to linearize $TT_{vt}$, Eqs.~\eqref{eq:eq_lin1}--\eqref{eq:eq_lin3} are added to the model incorporating the non-negative continuous variable $z_{ibjvt}$.

\begin{subequations}
	\begin{alignat}{2}
		&  z_{ibjvt} \leq  n_{bj}, {\begin{array}{l}
				\forall i,j \in I, b\in B, v\in V, t\in T
		\end{array}} \label{eq:eq_lin1}\\[0.05in]
		&  z_{ibjvt} \leq x_{ijvt},
		{\begin{array}{l} 
				\forall i,j \in I, b\in B, v\in V, t\in T
		\end{array}} \label{eq:eq_lin2}\\[0.05in]
		&  z_{ibjvt} \geq  n_{bj} + x_{ijvt} - 1,
		{\begin{array}{l}
				\forall i,j \in I, b\in B, v\in V, t\in T
		\end{array}} \label{eq:eq_lin3}
	\end{alignat}
\end{subequations}

Thus, $TT_{vt}$ can be rewritten as the linear expression of Eq.~\eqref{eq:eq_TT'}.

\begin{equation} 
	\label{eq:eq_TT'}
	TT'_{vt} = T_U \sum_{i\in {I}} x_{i0vt} + \sum_{i\in {I^0}} \sum_{j\in{I^0}} \left( C_{ij} x_{ijvt} + \sum_{b\in B} S_{b} z_{ibjvt} \right)
\end{equation}

When $TT_{vt}$ is replaced by $TT'_{vt}$ in the objective function~\eqref{eq:eq_2} and restriction~\eqref{eq:eq_2h} the model becomes linear. For showing the effectiveness of this linearization, both formulations (MIQP and MILP) were solved. Specifically, they were coded using pyomo~\citep{bynum2021pyomo} and solved with Gurobi v10.0.2~\citep{gurobi}.

\subsection{Genetic algorithm}
\label{sec:GA}	
Genetic Algorithm (GA) is a population-based metaheuristic method proposed by \cite{holland92} that is inspired by Charles Darwin's natural evolution and was designed to search the optimal solution in optimization problems. 

The basic elements of GA are the representation of the problem solutions and operators inspired by biological evolution: selection (of the individuals best adapted to the environment, i.e. to the problem), crossover or recombination and mutation. Selection is independent of representation, but the crossover and mutation operators are not. In general, GA moves between the set of solutions of the problem, called phenotype, and the set of individuals of a 'natural' population, encoding the information of each solution in a chain called chromosome. The symbols that form the string are called genes. When the representation in the chromosomes is made with strings of binary digits, it is known as genotype. Other types of representations are: real, integer, permutations, etc.

Evolution in GA usually starts from a population of individuals, or set of potential solutions to the problem, generated randomly, and continues with an iterative process, in which each new population is associated with a new generation. In each generation, the fitness of each individual in the population is evaluated and, depending on this value, candidate individuals are selected to be recombined, with a certain crossover rate, and whose descendants can be mutated or not, with a certain mutation rate, to give rise to a new generation. To avoid the loss of the best solution(s) found in a given generation, elitist strategies are usually applied that may or may not copy these elitist individuals in the next generation. Usually, the algorithm terminates when a maximum number of generations has been produced. 

\begin{algorithm}
	\caption{Mixed GA}\label{alg:GA}
	\begin{algorithmic}
		\Require problem data and algorithm parameters: population size ($N$), number of generations ($Ngen$), selection, crossover and mutation operators, crossover and mutation rates and number of elite solutions (individuals) considered ($\geq0$). 
		\State \textbf{Generate} a mixed initial population: $P^0 \gets \{(pop,mask)_n | 0 \le n < N \}$
		\State \textbf{Evaluate} the fitness function of each individual of the population.
		\State \textbf{Save} the elitist solutions.
		\State $ng = 0$
		\While{$ng < Ngen$}
		\State $ng = ng + 1$
		\State \textbf{Select} individuals to be crossed (2 times the population size).
		\For{$i$ in $range(0,N,2)$}
		\State \textbf{Cross} consecutive pairs (of both encodings) of the selected 
		\State individuals, with the fixed crossover rate.
		\State \textbf{Mute} the descendants (of both encodings) of the previous crossover 
		\State of two parents, with the fixed mutation rate.
		\State \textbf{Evaluate} the fitness function of the mutated descendants.
		\State \textbf{Insert} the mutated descendants into the new generation.	
		\EndFor	
		\State \textbf{Apply elitism} by adding the elitist solutions obtained in the previous 
		\State generation to the new generation.
		\State \textbf{Save} the elitist solutions of the new generation.
		\EndWhile
		\State \textbf{Save} the best individual found in a file.
	\end{algorithmic}
\end{algorithm}

Here, a mixed GA applied to a PCVRP is proposed to solve the problem described in Section~\ref{sec:problem_description} (see Algorithm~\ref{alg:GA}). Thus, two chromosomes are considered to encode the candidate solutions: $x \rightarrow (pop,mask)$ where $pop = [p_{it}]$ and $mask = [m_{it}]$ with $i\in I$ and $t\in T$, where $pop$ is a matrix whose rows correspond to the daily permutations of the order of MSW collection from the $n_I$ collection points, assuming that all of them are visited every day, and $mask$ is a binary matrix whose terms $m_{it}$ take the value 1, if on day $t \in (T-T')$ MSW is collected from collection point $i$, or 0 otherwise. This definition of $mask$ guarantees that restriction~\eqref{eq:eq_2e} is satisfied. Furthermore, by handling permutations, constraints \eqref{eq:eq_2d} and \eqref{eq:eq_2f} are always satisfied. Also, by not incorporating the depot in the permutations, it facilitates the constraint \eqref{eq:eq_2a} to be fulfilled. The remaining constraints are incorporated into the GA when the fitness function is evaluated. On the one hand, individuals are repaired by adding ones to the $mask$ when $t \in (T-T')$, in order to meet constraints \eqref{eq:eq_2b}, \eqref{eq:eq_2c} and \eqref{eq:eq_2i}~--~\eqref{eq:eq_2m}. On the other hand, they are penalized in case they do not meet any of the restrictions \eqref{eq:eq_2g} and~\eqref{eq:eq_2h}, i.e. if the number of daily routes exceeds the size of the vehicle fleet ($n_V$) and/or if a route cannot be completed within the drivers' working day, as shown in Eq.~\eqref{eq:eq_ff}.

\begin{equation} 
	\label{eq:eq_ff}
	f_{fitness}(x) = f(x) + \lambda \sum_{\substack{ t \in (T-T') \\ nr_t > n_V}} (nr_t - n_V)/n_V + \gamma \sum_{\substack{v \in V \\ t \in (T-T') \\ TT_{vt} > T_{L}}} (TT_{vt} - T_{L}) 
\end{equation}

{\noindent}where $f(x)$ is the objective function~\eqref{eq:eq_2}, $nr_t$ is the number of routes on day $t~\in~(T-T')$, $\lambda \in \mathbb{R}^{+}$ and $\gamma \in \mathbb{R}^{+}$.

During each successive generation, individuals are selected from the current population to produce a new generation. In this paper, a tournament selection that randomly samples two individuals without replacement and selects the one with the best fitness function value was considered and, for binary representation, two-point crossover and uniform mutation operators were chosen, as they have proven to be efficient operators in the literature. However, for the permutation representation, a design of experiment similar to the one performed by~\cite{Rossit2024} was performed.
 
As an illustrative example, the chromosome associated with the solution obtained for the $i.12.1$ instance with the MILP model is presented in the Appendix \ref{secA1}. 

\section{GA parametrization}\label{sec:GAparam}
This section reports the instances used in the paper as well as the $4\times3^3$ factorial design performed to tune the mixed GA proposed.

\subsection{Description of instances}\label{sec:instances}
\begin{figure}[]
	\centering
	\includegraphics{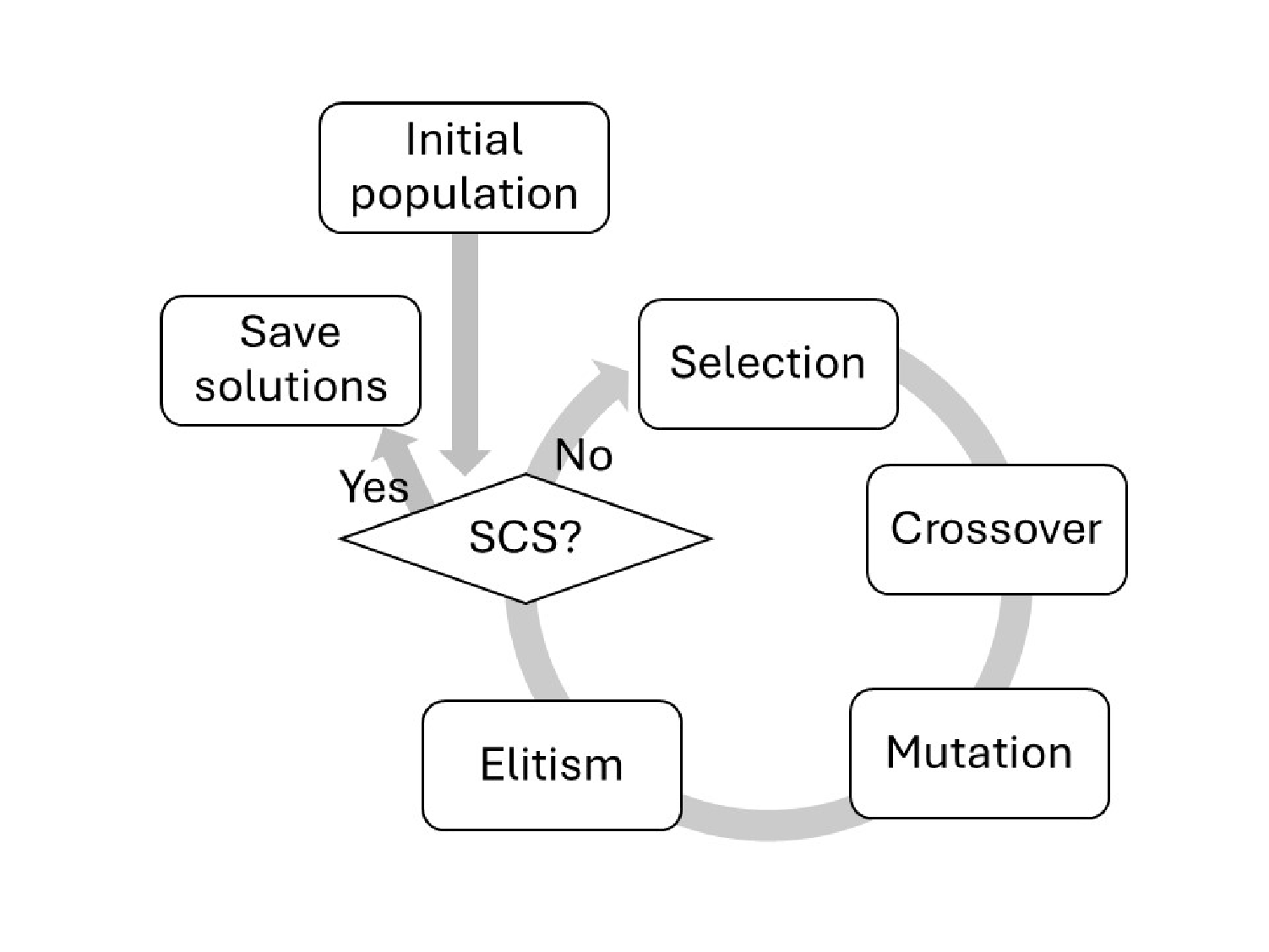}
	\caption{Collection points (in red) of instance $i.163.1$ (Bah\'ia Blanca, Argentina). The blue triangle marks the location of the depot.}\label{fig:fig1}
\end{figure}

This paper presents results obtained with several instances based on data collected in field studies in the city of Bah\'ia Blanca, Argentina~\citep{cavallin2020application}. The instances, which can be downloaded from Github\footnote{\url{https://github.com/diegorossit/ANOR-S-24-01950.git}}, have the generic name $i.n_I.m$, where $n_I$ represents the number of collection points and $m$ is a numbering of the instances in each particular case. Specifically, $m = 1,\dots, 5$ when $n_I = 12$ and $m = 1$ when $n_I \in \{40, 80, 120, 163\}$. For example, there is a single instance of 163 collection points whose name is $i.163.1$. Each instance consists of two files: \emph{time.txt}, which contains the matrix of travel times (minutes) from the depot to the collection points, between collection points, and from the collection points to the depot, and \emph{waste.txt}, which contains the geographic coordinates of each collection point, as well as the waste ($m^3$) generated daily by the urban area assigned to each collection point. Figure~\ref{fig:fig1} shows the 163 collection points of instance $i.163.1$. The smaller instances were randomly generated from this one. As an illustrative case, the daily routes obtained in the solution of the MILP model on instance $i.12.1$ are presented in the Appendix \ref{secA1}.

Table~\ref{tab:tableA3} shows the capacity, service time and weekly cost associated with the eight bin combinations considered in this study, which along with the specifications of the three types of bins that compose them, can be found in \cite{maheo2023solving}. However, unlike their approach, which assumes a uniform service time for all bin combinations, in this work the service time has been estimated based on the specific type of the bins constituting each combination. This modification adds more complexity to the objective problem and causes the appearance of bilinear terms in Eq.~\eqref{eq:eq_2}. The estimated lifespan of the bis is ten years~\citep{brogaard2012quantifying} and the maintenance cost of each bin is approximately 5\% of the purchase cost~\citep{d2016full}. To calculate the estimated weekly cost of each bin combination, the purchase and maintenance expenses for each bin were summed and divided by the total number of weeks of expected lifespan. Bins service times are estimated using the field work of~\cite{carlos2019influence}, who considered the collection vehicle fleet to be homogeneous. 

\begin{table}[h]
	\caption[]{Description of the bin combinations considered.}\label{tab:tableA3}%
	\begin{tabular}{@{}llll@{}}
		\toprule
		Bin & Capacity & Service time & Weekly cost \\
		combination & ($m^3$) & (minutes) & (US\$) \\
		\midrule
		0 & 1.1 & 0.70 & 0.78 \\ 
		1 & 2.2 & 1.40 & 1.56 \\ 
		2 & 2.4 & 0.66 & 1.56 \\ 
		3 & 3.3 & 2.10 & 1.56 \\ 
		4 & 3.5 & 1.36 & 1.56 \\ 
		5 & 4.3 & 1.37 & 1.56 \\ 
		6 & 4.8 & 1.32 & 1.56 \\ 
		7 & 5.6 & 1.33 & 1.56 \\ 
		\botrule
	\end{tabular}
\end{table}

As in~\cite{maheo2023solving}, since the instances are smaller than the actual collection zones in the city, vehicle capacity, fleet size and working day length are adjusted so that the problem does not become trivial, where a single vehicle can collect all the waste in a single trip. In this sense, a vehicle capacity of 12 $m^3$ has been considered for the instances of 12 collection points and 21 $m^3$ for the rest. Moreover, and in contrast to that work, this paper considers a longer planning horizon taking into account a whole week (seven days), with no waste collection scheduled on Sundays as this is a rest day for drivers. 

Finally, the estimated cost of the collection vehicles was set at US\$0.5764 per minute~\citep{d2016full}, and based on the paper of~\cite{Owusu-Nimo2019} a vehicle unloading time of 8 minutes was considered. In addition, since the mathematical model establishes that each vehicle can perform at most one route per day, the size of the fleet and the duration of the working day has been established by applying Eqs.~\eqref{eq:eq_nL} and~\eqref{eq:eq_TL}, respectively.

\begin{equation} 
	\label{eq:eq_nL}
	n_V = \lceil n_I/10 \rceil = \min \{k \in \mathbb{Z} | n_I/10 \leq k\}  
\end{equation}

\begin{equation} 
	\label{eq:eq_TL}
	T_L = \lceil \sum_{i,j \in I_0} C_{ij} / \left(n_V (n_V-1) |T-T'|\right) \rceil
\end{equation}

\subsection{Design of experiments}\label{sec:design}
A $4\times3^3$ factorial design on the $i.12.1$ instance was considered to tune the mixed GA, specifically focusing on the genetic operators associated to the representation with permutations. The objective of this design of experiments is to study the effect of four parameters of the GA on its performance and its runtime when searching for an optimal or quasi-optimal solution of the considered problem. The factors of this factorial design are: \textit{Crossover operator} (F1), \textit{Crossover rate} (F2), \textit{Mutation operator} (F3) and \textit{Mutation rate} (F4). Their respective levels are:

\begin{itemize}
	\item [] \textbf{F1}: (1) Partially Mapped Crossover (PMX)~\citep{goldberg14}, (2) Order Crossover (OX)~\citep{davis85}, (3) Cicle Crossover (CX)~\citep{oliver87} and (4) Modified Cicle Crossover (CX2)~\citep{hussain17}.
	\item [] \textbf{F2}: Crossover rate: (1) 0.8, (2) 0.85 and (3) 0.9.
	\item [] \textbf{F3}: Mutation operator: (1) Exchange mutation (EM), (2) Insertion mutation (IM) and (3) Inversion mutation (INM)~\citep{Rossit2024}.
	\item [] \textbf{F4}: Mutation rate: (1) 0.05, (2) 0.10 and (3) 0.15.
\end{itemize}

In addition, populations of 100 individuals with elitism of 2 and 1000 generations were fixed. For binary encoding, the crossover rate is the same as permutation encoding and the mutation rate is equal to 1/number of collection points. Finally, so that the penalties for infeasible solutions in the fitness function defined by Eq.~\eqref{eq:eq_ff} would not be too high, we set $\lambda = 100 $ and $\gamma = 1000$. 

For each combination of factor levels, referred to as treatment, 30 independent runs were carried out. A total of 3,240 runs of the mixed GA were performed. Runs were performed on an Intel(R) Core(TM) i7-10700K CPU @ 3.80GHz 3.79 GHz with 16.0 GB of RAM.

\begin{figure}[]
	\centering
	\includegraphics[scale=0.65]{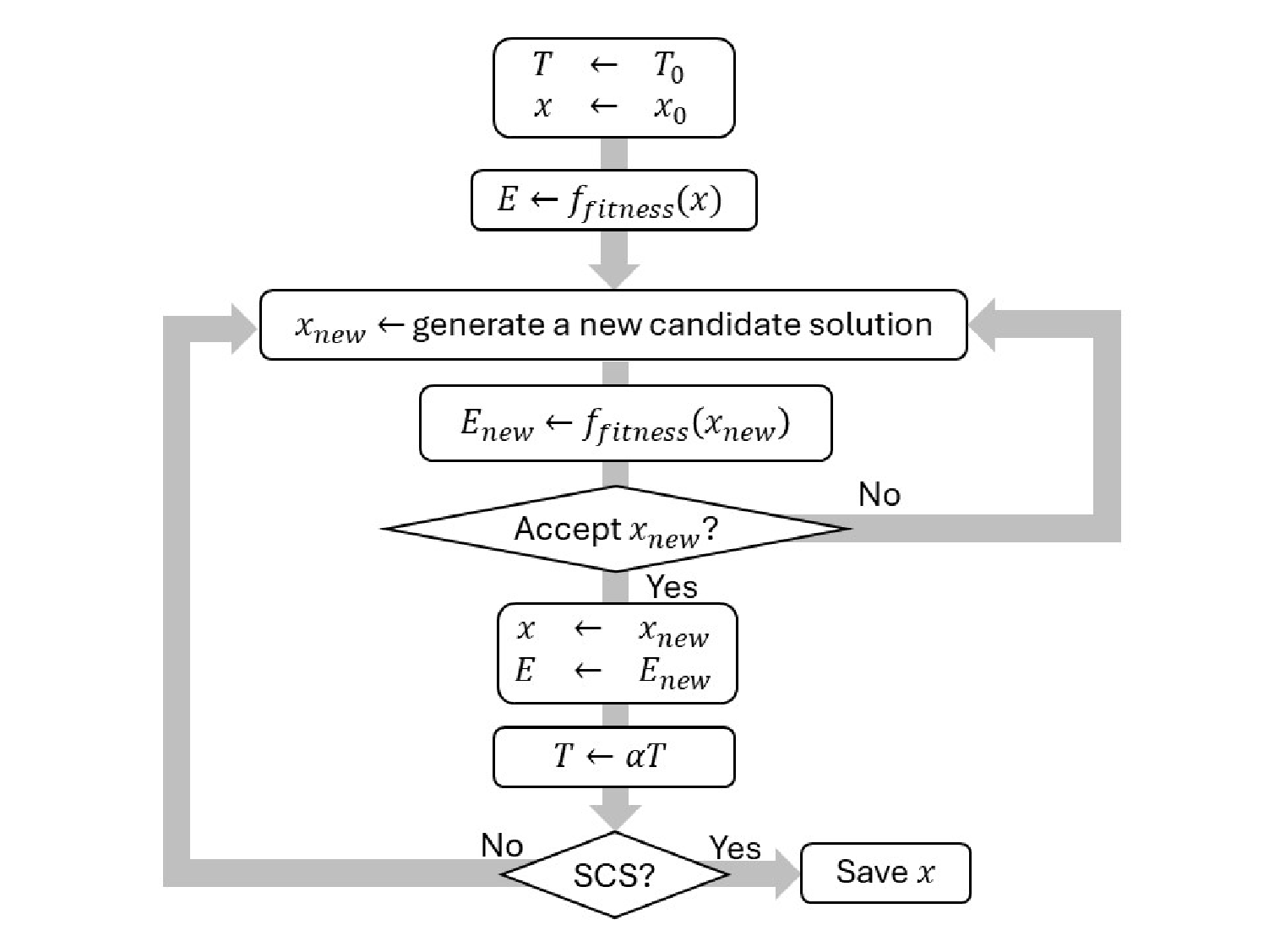}
	\caption{Main effects when the response variable is the overall cost (US\$).}\label{fig:fig2}
\end{figure}

In the context of the design of experiments, main effects are defined as the changes in the mean of the response variable that are due to the individual action of each factor. From Figure~\ref{fig:fig2} it could be inferred that, on overage, the best treatment would be: (F1~=~CX, F2~=~0.8, F3~=~EM, F4~=~0.05). However, this conclusion needs to be confirmed with the corresponding hypothesis tests.

A Kruskal-Wallis Rank Sum Test~\citep{hollander73} has been carried out to establish whether there are significant differences between the levels of each factor. When the value of a Kruskal-Wallis test is significant ($p$-value $< 0.05$), a multiple comparison test after Kruskal-Wallis~\citep{siegel88} between treatments has been carried out to identify which levels are different with pairwise comparisons adjusted appropriately for multiple comparisons. The tests were performed in RStudio Desktop~\citep{RStudio2023}, a coding environment for R~\citep{RCoreTeam2022}, a free software environment for statistical computing and graphics. Table~\ref{tab:table1} shows the $p$-values obtained when the \textbf{kruskal.test()} function of the R package \textit{stats}~\citep{RCoreTeam2022} is run on the overall costs obtained. It can be concluded that, for each factor, there are at least a pair of levels that are significantly different. Table~\ref{tab:table1} also shows the output of the multiple comparison test after Kruskal-Wallis when the \textbf{kruskalmc()} function of the R package \textit{pgirmess}~\citep{giraudoux23} is run on the overall costs obtained: if TRUE (T), then statistically significant differences are found between the compared levels. Therefore, it can be seen that this test detected significant differences between all pairs of levels of factors F1 and F2, but only detected significant differences between pair 1--3 of factor F3 and pairs 1--2 and 1--3 of factor F4. That is, hypothesis tests support that the best performing treatment is: (F1~=~CX, F2~=~0.8, F3~=~EM, F4~=~0.05), but also (F1~=~CX, F2~=~0.8, F3~=~IM, F4~=~0.05). 

\begin{table}[]
	\caption{Overall cost: summary of statistical tests}\label{tab:table1}%
	\begin{tabular}{@{}lllll@{}}
		\toprule
		\multicolumn{5}{c}{Kruskal-Wallis rank sum test: $p$-value} \\ 
		\midrule
		\multicolumn{2}{l}{F1: Crossover operator} & F2: Crossover rate & F3: Mutation operator & F4: Mutation rate \\
		\midrule
		\multicolumn{2}{l}{$< 2.2e-16$} & $< 2.2e-16$ & 0.0004956 & 0.0005319 \\ 
		\midrule
		\multicolumn{5}{c}{The multiple comparison test after Kruskal-Wallis} \\ 
		\midrule
		\multicolumn{2}{l}{F1: Crossover operator} & F2: Crossover rate & F3: Mutation operator & F4: Mutation rate \\
		\midrule
		1--2 T & \qquad 2--3 T & 1--2 T & 1--2 F & 1--2 T \\ 
		1--3 T & \qquad 2--4 T & 1--3 T & 1--3 T & 1--3 T \\
		1--4 T & \qquad 3--4 T & 2--3 T & 2--3 F & 2--3 F \\
		\botrule
	\end{tabular}
\end{table}

\begin{figure}[]
	\centering
	\includegraphics[scale=0.65]{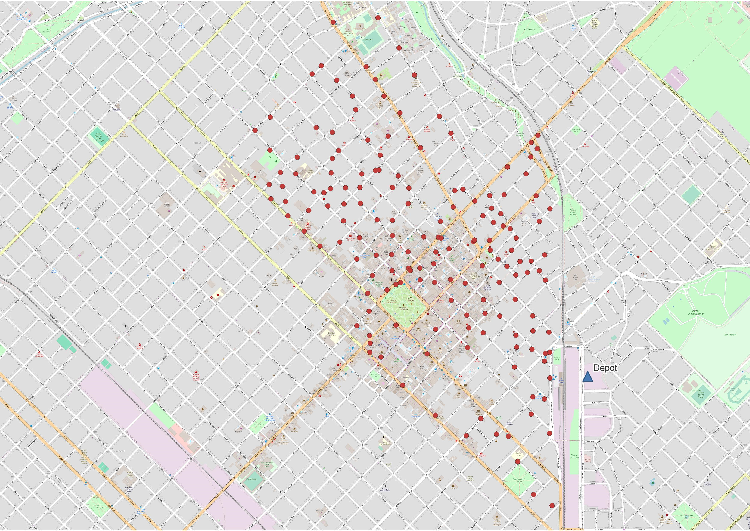}
	\caption{Main effects when the response variable is the mixed GA runtime (seconds).}\label{fig:fig3}
\end{figure}

Main effects, when the response variable is the mixed GA runtime (seconds), is shown in Figure~\ref{fig:fig3}. Table~\ref{tab:table2} shows the $p$-values obtained when performing a Kruskal-Wallis rank sum test on the  mixed GA runtime. It can be concluded that, for each factor, there are at least a pair of levels that are significantly different. Table~\ref{tab:table2} also shows that the multiple comparison test after Kruskal-Wallis detected significant differences between all pairs of levels of factors F1, F2 and F4, but only detected significant differences between pair 1--3 of factor F3. Consequently, the hypothesis tests support that the crossover operator that needs the least runtime is CX2 closely followed by CX (see Figure~\ref{fig:fig3}). However, the performance of CX is significantly better than that of CX2 (see Figure~\ref{fig:fig2}).

\begin{table}[]
	\caption{Mixed GA runtime: summary of statistical tests}\label{tab:table2}%
	\begin{tabular}{@{}lllll@{}}
		\toprule
		\multicolumn{5}{c}{Kruskal-Wallis rank sum test: $p$-value} \\ 
		\midrule
		\multicolumn{2}{l}{F1: Crossover operator} & F2: Crossover rate & F3: Mutation operator & F4: Mutation rate \\
		\midrule
		\multicolumn{2}{l}{$< 2.2e-16$} & $< 2.2e-16$ & 0.0007906 & $< 2.2e-16$ \\		
		\midrule
		\multicolumn{5}{c}{The multiple comparison test after Kruskal-Wallis} \\ 
		\midrule
		\multicolumn{2}{l}{F1: Crossover operator} & F2: Crossover rate & F3: Mutation operator & F4: Mutation rate \\
		\midrule
		1--2 T & \qquad 2--3 T & 1--2 T & 1--2 F & 1--2 T \\
		1--3 T & \qquad 2--4 T & 1--3 T & 1--3 T & 1--3 T \\
		1--4 T & \qquad 3--4 T & 2--3 T & 2--3 F & 2--3 T \\
		\botrule
	\end{tabular}
\end{table}

\section{Results and discussion}
\label{sec:NumResults}
This section presents the results obtained both for comparing the exact models (MIQP and MILP) with the mixed GA, with the five instances of 12 collection points, and the solutions obtained by running the mixed GA with the remaining instances considered. In view of the conclusions made in the previous section, all mixed GA runs were performed with treatment: (F1~=~CX, F2~=~0.8, F3~=~EM, F4~=~0.05).

\begin{table}[]
	\caption{Instances $i.12.*$: Statistical summary of the overall costs over 30 independent mixed GA runs.}\label{tab:GA_results}%
	\begin{tabular}{@{}llllllll@{}}
		\toprule
		Instance & Min & Median & Mean & \multicolumn{2}{l}{Confidence interval} & Runtime (s)& S-W test \\
		\midrule
		$i.12.1$ & 194.6 & 205.1 & 205.1 & 204.2 & 206.3 & 922.9 & 0.0348  \\
		$i.12.2$ & 198.5 & 203.8 & 204.2 & 202.9 & 205.5 & 927.0 & 0.1724  \\
		$i.12.3$ & 196.5 & 206.2 & 206.2 & 205.1 & 207.1 & 910.2 & 0.0462  \\
		$i.12.4$ & 187.5 & 195.8 & 195.6 & 194.4 & 196.8 & 903.4 & 0.1254  \\ 
		$i.12.5$ & 192.8 & 196.0 & 196.7 & 195.8 & 197.6 & 917.4 & 0.3189  \\ 
		\botrule
	\end{tabular}
\end{table}

Table~\ref{tab:GA_results} summarizes the statistical results obtained from a sample of 30 independent GA runs for the five instances of 12 collection points. The parameters considered were the same as those for GA tuning, but this time 10,000 generations were performed so that, on the one hand, the GA convergence curve would stabilize and, on the other hand, the execution time of each instance (30 independent mixed GA runs) would be of the order of the execution time of the exact methods (8 hours). From left to right, Table~\ref{tab:GA_results} shows, for each sample of the 30 best overall costs (US\$) obtained with the 30 independent mixed GA runs: the minimum, median, mean, 95 percent confidence interval of the mean or of the (pseudo)median, as well as the mean run time (in seconds) and the $p$-value obtained in the Shapiro Wilk (S-W) test with the \textbf{shapiro.test()} function. If a distribution is symmetric, then the pseudomedian and median coincide~\citep{hollander73}. When the value of a S-W test is significant ($p$-value $< 0.05$) then the distribution of the sample is significantly different from a normal distribution. Based on this $p$-value, a 95 percent confidence interval is calculated for the (pseudo)median with the \textbf{wilcox.test()} function, if $p$-value $< 0.05$, or for the mean with the \textbf{t.test()} function, otherwise. These three test functions can be found in the R package \textit{stats}~\citep{RCoreTeam2022}.

\begin{table}[]
	\caption{Instances $i.12.*$: MILP and MIQP solutions.}\label{tab:MILP_results}%
	\begin{tabular}{@{}llllll@{}}
		\toprule
		\multicolumn{6}{c}{MILP} \\
		\midrule
		Instance & Bin cost  & Routing cost & Overall cost & Lower bound & Optimality gap \\
		\midrule
		$i.12.1$ & 45.47 & 157.35 & 202.82 & 141.66 & 30.15\% \\ 
		$i.12.2$ & 39.64 & 150.11 & 189.75 & 131.41 & 30.75\% \\ 
		$i.12.3$ & 43.25 & 153.40 & 196.65 & 138.14 & 29.75\% \\ 
		$i.12.4$ & 42.94 & 149.46 & 192.40 & 136.50 & 29.06\% \\ 
		$i.12.5$ & 42.85 & 146.73 & 189.58 & 128.63 & 32.15\% \\ 
		\toprule
		\multicolumn{6}{c}{MIQP} \\
		\midrule
		Instance & Bin cost  & Routing cost & Overall cost & Lower bound & Optimality gap \\
		\midrule
		$i.12.1$ & 40.35 & 148.69 & 188.97 & 139.56 & 26.15\% \\ 
		$i.12.2$ & 39.18 & 152.03 & 191.21 & 129.38 & 32.34\% \\ 
		$i.12.3$ & 48.74 & 153.36 & 202.10 & 140.40 & 30.53\% \\ 
		$i.12.4$ & 40.26 & 144.75 & 185.01 & 133.71 & 27.73\% \\ 
		$i.12.5$ & 41.80 & 145.11 & 186.91 & 127.96 & 31.54\% \\ 
		\botrule
	\end{tabular}
\end{table}

Table~\ref{tab:MILP_results} shows the results obtained when solving the instances of 12 collection points with the MILP and MIQP models. From left to right, this table shows the bin installation and maintenance cost, the MSW collection cost, the overall cost, the lower bound estimated by Gurobi (dual objective bound) and the optimality gap reported by Gurobi, which is calculated with Eq.~\eqref{eq:gurobi_gap}.

\begin{equation}
	\label{eq:gurobi_gap}
	\text{Optimality gap} = \frac{|\text{Overall cost} - \text{lower bound}|}{|\text{Overall cost}|}
\end{equation}

The optimality gap is useful in computationally complex problem, such as the optimization problem addressed in this paper, where it is usual for exact solvers to report feasible but not proven optimal solutions. In these cases, this gap, which is estimated with the lower bound, indicates an estimate of how far the obtained feasible solution might be from the optimal solution. Further details on how Gurobi computes the optimality gap can be found at~\cite{gurobi}.

\begin{table}[]
	\caption{Instances $i.12.*$: GA versus MILP and MIQP}\label{tab:GA_versus_MILP}%
	\begin{tabular}{@{}lllllllll@{}}
		\toprule
		&  & Cost & $i.12.1$ & $i.12.2$ & $i.12.3$ & $i.12.4$ & $i.12.5$ & Mean \\ 
		\midrule
		\multicolumn{1}{c}{\multirow{3}{*}{}} & \multirow{3}{*}{GA Min}  & Bin ($C_B$) & 51.95 & 47.26 & 51.95 & 51.61 & 53.83 \\
		\multicolumn{1}{c}{}	&	& Routing ($C_R$) & 142.65 & 151.25 & 144.54 & 135.90 & 138.94 \\
		\multicolumn{1}{c}{}	&	& Overall ($C_O$) & 194.60 & 198.51 & 196.49 & 187.51 & 192.77 \\
		\midrule
		\multirow{6}{*}{\rotatebox{90}{MILP vs GA}}  & \multirow{3}{*}{GA Min}  & \% dif. $C_B$ & 14.25 & 19.22 & 20.12 & 20.19 & 25.62 & 19.88   \\
		&  & \% dif. $C_R$ & -9.34 & 0.76 & -5.78 & -9.07 & -5.31& -5.75   \\
		&  & \% dif. $C_O$ & -4.05 & 4.62 & -0.08 & -2.54 & 1.68 & -0.08  \\
		\cmidrule{2-9}
		& GA & \% dif. $C_B$ & 13.57 & 21.69 & 21.90 & 18.95 & 17.16 & 18.65\\
		& Mean/ & \% dif. $C_R$ & -2.17 & 3.80 & 0.20 & -3.35 & -0.29 & -0.36   \\
		& Median & \% dif. $C_O$  & 1.08 & 7.53 & 4.86 & 1.63 & 3.65 & 3.75   \\
		\midrule
		\multirow{6}{*}{\rotatebox{90}{MIQP vs GA}}  & \multirow{3}{*}{GA Min}  & \% dif. $C_B$ & 28.75 & 20.62 & 6.59 & 28.19 & 28.78 & 22.59    \\
		&  & \% dif. $C_R$ & -4.06 & -0.51 & -4.93 & -10.61 & -8.61 & -5.74   \\
		&  & \% dif. $C_O$ & 2.98 & 3.82 & -2.78 & 1.35 & 3.14& 1.70   \\
		\cmidrule{2-9}
		& GA & \% dif. $C_B$ & 27.98 & 23.12 & 8.17 & 26.86 & 20.10 & 21.24 \\
		& Mean/ & \% dif. $C_R$ & 3.52 & 2.48 & 0.22 & -0.20 & 0.82 & 1.37   \\
		&  Median &  \% dif. $C_O$  & 8.49 & 6.71 & 2.03 & 5.69 & 5.13 & 5.61   \\
		\botrule
	\end{tabular}
\end{table}

The comparison of the results of Gurobi with both the best solution obtained with the mixed GA and the corresponding mean or median, depending on whether or not there is sufficient evidence to assume that the sample is normally distributed (see the $p$-value of the Shapiro-Wilk test in Table~\ref{tab:GA_results}), are shown in Table~\ref{tab:GA_versus_MILP}, where the percentage difference between the costs obtained with the exact methods and the GA has been calculated with Eq.~\eqref{eq:dif}:

\begin{equation}
	\label{eq:dif}
	\text{\% dif. } C_{*} = (C_{*}^{GA} - C_{*}^{MILP/MIQP}) / C_{*}^{MILP/MIQP})
\end{equation}

{\noindent}where $C_{*}^{GA}$ is the cost obtained with the GA and $C_{*}^{MILP/MIQP}$ is the cost obtained with the MILP or MIQP model.

When comparing the best solution of the GA with the MILP model, GA consistently obtains a better routing cost and an average improvement of 5.75\%, getting the biggest difference in the $i.12.1$ instance (-9.34\%). However, GA falls short in obtaining low bin cost: the MILP model is able to obtain much better results, with an average improvement of 19.88\% (and up to 25.62\% in instance $i.12.5$). Regarding the overall cost, the best GA solution outperforms the MILP model in three of the five instances. A similar behaviour is observed when comparing the GA with the MIQP model. The best GA solution is really efficient in constructing routes, outperforming the MIQP model in all instances. The mean improvement is 5.74\% with up to 10.61\% in instance $i.12.4$. The MIQP solution is much better in the bin cost, obtaining up to 28.78\% smaller cost that the GA. Overall, the MIQP is usually less expensive (an average of 1.70\%) except in instance $i.12.3$, in which GA excels MIQP. 

When comparing the MILP and MIQP solutions with the mean/median cost values of the AG solutions, the situation is similar. For example, GA obtains, in general, better solution in terms of routing cost than the MILP model (0.36\% better on average), but worse solutions in terms of bin costs (18.65\% worse on average ). However, the improvement in routing cost is relatively small and, thus, the comparison of the mean overall cost of the GA solutions and the MILP solution is on average 3.75\% worse. When comparing the mean values of the GA solutions to the MIQP solucion, it is almost always better in every cost for all instances except for the routing cost in the $i.12.4$ instance (-0.20\%).

In any case, the proposed mixed GA has proven to be competitive in solving computationally complex problems, such as the problem modeled in this paper, since without having allowed it to evolve beyond 10,000 generations with 100 individuals, the average percentage difference between the average overall costs obtained with the GA and the MILP or MIQP solution is, respectively, less than 4\% in the case of MILP and less than 6\% in the case of MIQP. 

Table~\ref{tab:GA_big_instances} summarizes the statistical results obtained from a sample of 30 independent GA runs for the instances of 40, 80, 120 and 163 collection points. The parameters considered were the same as those for GA tuning, but this time 10,000 generations were performed in the case of 40 collection points, and 20,000 generations in the other cases, so that the GA convergence curve would stabilize. Furthermore, since the value of the objective function is higher the larger the instance considered, different values of $\lambda$ were considered to evaluate the fitness function \eqref{eq:eq_ff}. These values are reported in the table. Specifically, Table~\ref{tab:GA_big_instances} shows, for each sample of the 30 best overall costs (US\$) obtained with the 30 independent mixed GA runs: the minimum, median, mean, 95 percent confidence interval of the mean or of the (pseudo)median, as well as the mean runtime (in hours) and the $p$-value obtained in the Shapiro Wilk (S-W) test with the \textbf{shapiro.test()} function. Again, as in Table~\ref{tab:GA_results}, based on this $p$-value, the 95 percent confidence interval is calculated for the (pseudo)median with the \textbf{wilcox.test()} function, if $p$-value $< 0.05$, or for the mean with the \textbf{t.test()} function, otherwise. It should be noted that feasible solutions were obtained in all runs and a linear growth of the average runtimes is observed as the number of collection points increases, with a slope around 0.035.

\begin{table}[]
	\caption{Mixed GA. Statistical summary of the overall costs obtained with the instances of 40, 80, 120 and 163 collection points. In each case, 30 independent runs have been performed. The average runtime is in hours.}\label{tab:GA_big_instances}%
	\begin{tabular}{@{}lllllllll@{}}
		\toprule 
		Instance & Min & Median & Mean & \multicolumn{2}{l}{Confidence interval} &  Runtime (h) & S-W test \\
		\midrule
		$i.40.1 (\lambda = 500)$  & 527 & 543 & 543 & 540 & 545 & 0.72 & 0.0080  \\
		$i.80.1 (\lambda = 1000)$  & 1077 & 1106 & 1106 & 1100 & 1111 & 2.84 & 0.9876  \\
		$i.120.1 (\lambda = 5000)$ & 1658 & 1726 & 1730 & 1720 & 1740 & 4.09 & 0.6069  \\
		$i.163.1 (\lambda = 10000)$ & 2358 & 2390 & 2395 & 2386 & 2403 & 5.58 & 0.2697 \\ 
		\botrule
	\end{tabular}
\end{table}

\section{Conclusion and future work}\label{sec:conclusions}

Effective MSW systems are essential for contemporary cities, playing a critical role in enhancing the livability and sustainability of smart city initiatives. However, careful planning and operation are paramount to achieve this goal. This paper addresses, in an integrated manner, two fundamental and interrelated problems in waste management: the sizing of collection points and the design of weekly collection route. These problems are often treated separately due to their inherent computational complexity. 

In this work, two mathematical formulations for this problem are presented: a Mixed-Integer Quadratic model and a Mixed-Integer Linear model. In addition, a genetic algorithm with a complex chromosome representation is proposed to handle the integrated problem. Regarding GA, extensive computational experiments were conducted, evaluating different crossover and mutation operators with various probabilities to identify optimal configurations for the target problem. Subsequently, a comparison was performed between the mathematical formulations and the GA, using small realistic instances, which proved the competitiveness of the proposed GA to solve the considered problem, since it obtained results not very far from those obtained with the exact formulations, with a limited number of evaluations of the objective function so that the execution time of the 30 runs of the GA, with each instance, was of the order of that assigned to the exact formulations (8 hours). Finally, the proposed GA was used to address larger real-world instances, obtaining feasible solutions in all runs and thus demonstrating its efficiency in tackling practical instances where mathematical models become intractable.Furthermore, a linear growth in the average runtimes is observed as the number of collection points increases, with a slope around 0.035. 

Future research endeavors include enhancing the mathematical formulations with valid cuts to improve performance, exploring other metaheuristics such as simulated annealing for comparative analysis, and potentially addressing the problem as a bi-objective optimization, to separately optimize installation and routing costs. The problem of increased runtimes when using metaheuristics can be tackled from the point of view of parallelization of the algorithms.

\bmhead{Acknowledgements}
The authors acknowledge funding by Consejer\'ia de Econom\'ia, Industria, Comercio y Conocimiento, Gobierno de Canarias (ES). Project code: CEI2021-05 (ULPGC).

%The authors are grateful for the financial support granted by the Consejer\'ia de Econom\'ia, Industria, Comercio y Conocimiento of the Government of the Canary Islands through the direct grant awarded to the ULPGC called "Support for R+D+i activity. Campus of International Excellence CEI CANARIAS-ULPGC", 319RT0574 - Red Iberoamericana Industria 4.0 of CYTED through research projects, the PICT-2021-I-INVI00217 of the Agencia I+D+i of Argentina and the PIBAA 0466CO of CONICET.Un

\section*{Declarations}
\bmhead{Funding} Open access funding provided by Universidad de Las Palmas de Gran Canaria.

\bmhead{Conflict of interest} The authors declare that they have no conflicts of interest.

\bmhead{Materials availability} This paper presents results obtained with several instances based on data collected in field studies in the city of Bah\'ia Blanca, Argentina~\citep{cavallin2020application}, and they can be downloaded from \url{https://github.com/diegorossit/ANOR-S-24-01950.git}.

%\bmhead{Author contribution} 

\begin{appendices}

\section{An illustrative example of a feasible solution for the waste collection problem}\label{secA1}
\begin{figure}[]
	\centering
	\includegraphics{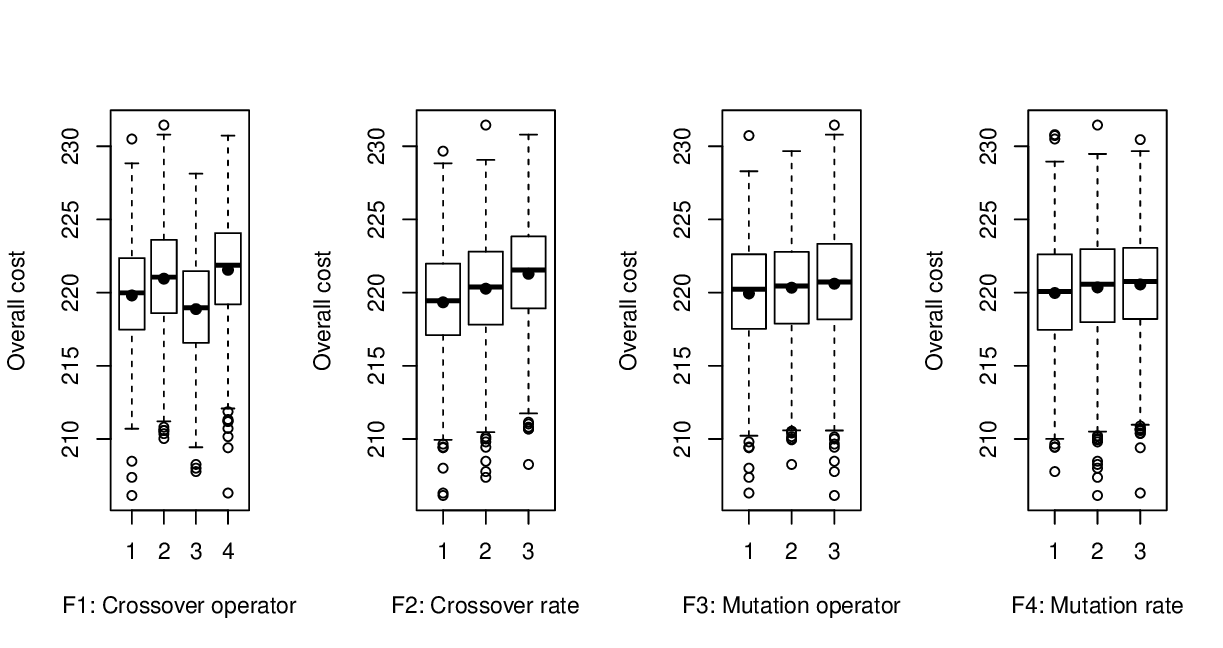}
	\caption{Waste collection routes associated with the solution shown in Table~\ref{tab:tableA1}. The time taken to complete them (in minutes) is also shown.}\label{fig:fig4}
\end{figure}

This section shows an illustrative (feasible) solution of the $i.12.1$ to analyze how the decodification function of the Genetic Algorithm works. Table~\ref{tab:tableA1} presents the chromosome of a feasible solution for this instance of twelve collection points and a planning horizon of one week in which the drivers' day of rest is Sunday. Firstly, the selection of the bin combination for a given collection point $i$ is made taking into account the maximum daily amount of waste accumulated in the planning horizon $w_i^{max} = \max_{t \in T} w_{it}$. For example, $w_1^{max} = 5.08$, therefore, at collection point 1 it will be necessary to install a bin combination with a capacity of 5.08 $m^{3}$ or more, in this case, the one identified as 7 (see Table~\ref{tab:tableA3}). For establishing the routes performed each day $t$ we follow the order of visitation of the collection points set by the chromosome $[y_{it}]$ and taking into account that $m_{it} = 0$ indicates that the collection point $i$ is not visited on day $t$, therefore, this point must be ignored in the order of visitation. For example, on Monday only collection points 7, 6, 12, 10, 3, 2 and 9 will be visited, in that order. The number of routes for these visits is established taking into account the capacity of the collection vehicles (12 $m^{3}$ in this case). A route ends when the vehicle does not have the capacity to empty the bin combination of the next collection point. Fig.~\ref{fig:fig4} shows the routes associated with the solution shown in Table~\ref{tab:tableA1} and Table~\ref{tab:tableA2} shows the waste collected by the vehicle on each route of that solution.

\begin{table}[]
	\caption[]{A feasible solution (chromosome) solution of instance $i.12.1$.}\label{tab:tableA1}%
	\begin{tabular}{@{}llllllllllllll@{}}
		\toprule
		Collection  & Bin & \multicolumn{3}{c}{MON} & & \multicolumn{3}{c}{TUE} & &  \multicolumn{3}{c}{WED} \\
		\cmidrule{3-5} \cmidrule{7-9} \cmidrule{11-13}
		points & combination & $p_{i0}$ & $m_{i0}$ & $w_{i0}$  & &  $p_{i1}$ & $m_{i1}$ & $w_{i1}$  & &  $p_{i2}$ & $m_{i2}$  & $w_{i2}$  \\
		\midrule
		1         & 7    & 8    & 0      & 2.54    & & 0    & 0      & 3.81     & &  4    & 1      & 5.08    \\
		2         & 7    & 6    & 1      & 3.24    & & 4    & 1      & 1.62     & &  0    & 0      & 1.62    \\
		3         & 2    & 5    & 1      & 2.34    & & 3    & 1      & 1.17     & &  2    & 1      & 1.17    \\
		4         & 6    & 9    & 0      & 2.98    & & 2    & 1      & 4.47     & &  5    & 0      & 1.49    \\
		5         & 6    & 10   & 0      & 3.18    & & 5    & 1      & 4.77     & &  6    & 0      & 1.59    \\
		6         & 5    & 2    & 1      & 2.42    & & 8    & 0      & 1.21     & &  7    & 0      & 2.42    \\
		7         & 7    & 1    & 1      & 5.28    & & 9    & 0      & 1.32     & &  8    & 0      & 2.64    \\
		8         & 7    & 11   & 0      & 3.69    & & 6    & 1      & 4.92     & &  9    & 0      & 1.23    \\
		9         & 4    & 7    & 1      & 3.16    & & 10   & 0      & 1.58     & &  3    & 1      & 3.16    \\
		10        & 2    & 4    & 1      & 2.34    & & 11   & 0      & 1.17     & &  1    & 1      & 2.34    \\
		11        & 5    & 0    & 0      & 2.00    & & 1    & 1      & 3.00     & &  10   & 0      & 1.00    \\
		12        & 4    & 3    & 1      & 2.66    & & 7    & 1      & 1.33     & &  11   & 0      & 1.33	\\ 
		\midrule
		Collection  & Bin & \multicolumn{3}{c}{THUR} & & \multicolumn{3}{c}{FRI} & & \multicolumn{3}{c}{SAT} \\
		\cmidrule{3-5} \cmidrule{7-9} \cmidrule{11-13}
		points & combination & $p_{i3}$ & $m_{i3}$ & $w_{i3}$  & & $p_{i4}$ & $m_{i4}$ & $w_{i4}$ & &  $p_{i5}$ & $m_{i5}$  & $w_{i5}$  \\
		\midrule
		1         & 7    & 5    & 0      & 1.27    & & 7    & 0      & 2.54     & &  4    & 1      & 3.81    \\
		2         & 7    & 6    & 0      & 3.24    & & 3    & 1      & 4.86     & &  3    & 1      & 1.62    \\
		3         & 2    & 7    & 0      & 1.17    & & 2    & 1      & 2.34     & &  2    & 1      & 1.17    \\
		4         & 6    & 8    & 0      & 2.98    & & 1    & 1      & 4.47     & &  8    & 1      & 1.49    \\
		5         & 6    & 9    & 0      & 3.18    & & 5    & 1      & 4.77     & &  9    & 1      & 1.59    \\
		6         & 5    & 3    & 1      & 3.63    & & 8    & 0      & 1.21     & &  6    & 1      & 2.42    \\
		7         & 7    & 2    & 1      & 3.96    & & 9    & 0      & 1.32     & &  11   & 0      & 2.64    \\
		8         & 7    & 10   & 0      & 2.46    & & 6    & 1      & 3.69     & &  0    & 0      & 1.23    \\
		9         & 4    & 11   & 0      & 1.58    & & 4    & 1      & 3.16     & &  10   & 1      & 1.58    \\
		10        & 2    & 1    & 1      & 1.17    & & 10   & 0      & 1.17     & &  1    & 1      & 2.34   \\
		11        & 5    & 0    & 0      & 2.00    & & 11   & 0      & 3.00     & &  7    & 1      & 4.00    \\
		12        & 4    & 4    & 1      & 2.66    & & 0    & 0      & 1.33     & &  5    & 1      & 2.66    \\
		\botrule
	\end{tabular}
\end{table}

\begin{table}[]
	\caption[]{Time and waste collected by the vehicle on each route of the solution shown in Table~\ref{tab:tableA1}. Index $g_i$ represents the previous collection point on the route (or the depot, if $i$ is the first point on the route). The collection vehicle capacity is 12 $m^{3}$.}\label{tab:tableA2}%
	\begin{tabular}{@{}lllllllllll@{}}
		\toprule
		\multirow{6}{*}{MON}  & \multirow{3}{*}{R1} & $i$     & 0 & 7    & 6    & 12   & 0     &   &    & \multirow{3}{*}{25.04 min} \\
		\cmidrule{3-8}
		&                     & $w_{i0}$  &   & 5.28 & 2.42 & 2.66 &       &       &     &  \\
		&                     & $v_{g_{i}i10}$ &   & 0.00 & 5.28 & 7.70 & 10.36 &       &     &  \\
		\cmidrule{2-11}
		& \multirow{3}{*}{R2} & $i$    & 0 & 10   & 3    & 2    & 9     & 0     &     & \multirow{3}{*}{23.05 min}  \\
		\cmidrule{3-9}
		&                     & $w_{i0}$   &   & 2.34 & 2.34 & 3.24 & 3.16  &       &    &   \\
		&                     & $v_{g_{i}i20}$ &   & 0.00 & 2.34 & 4.68 & 7.92  & 11.08 &      & \\
		\midrule
		\multirow{6}{*}{TUE}  & \multirow{3}{*}{R1} & $i$     & 0 & 11   & 4    & 3    & 2     & 0     &    & \multirow{3}{*}{26.00 min}   \\
		\cmidrule{3-9}
		&                     & $w_{i1}$   &   & 3.00 & 4.47 & 1.17 & 1.62  &       &     &  \\
		&                     & $v_{g_{i}i11}$ &   & 0.00 & 3.00 & 7.47 & 8.64  & 10.26 &     &  \\
		\cmidrule{2-11}
		& \multirow{3}{*}{R2} & $i$     & 0 & 5    & 8    & 12   & 0     &       &    & \multirow{3}{*}{22.29 min}   \\
		\cmidrule{3-8}
		&                     & $w_{i1}$   &   & 4.77 & 4.92 & 1.33 &       &       &     &  \\
		&                     & $v_{g_{i}i21}$ &   & 0.00 & 4.77 & 9.69 & 11.02 &       &      & \\
		\midrule
		\multirow{3}{*}{WED}  & \multirow{3}{*}{R1} & $i$     & 0 & 10   & 3    & 9    & 1     & 0     &    & \multirow{3}{*}{25.80 min}   \\
		\cmidrule{3-9}
		&                     & $w_{i2}$   &   & 2.34 & 1.17 & 3.16 & 5.08  &       &     &  \\
		&                     & $v_{g_{i}i12}$ &   & 0.00 & 2.34 & 3.51 & 6.67  & 11.75 &    &   \\
		\midrule
		\multirow{3}{*}{THUR} & \multirow{3}{*}{R1} & $i$     & 0 & 10   & 7    & 6    & 12    & 0     &   & \multirow{3}{*}{26.00 min}    \\
		\cmidrule{3-9}
		&                     & $w_{i3}$   &   & 1.17 & 3.96 & 3.63 & 2.66  &       &     &  \\
		&                     & $v_{g_{i}i13}$ &   & 0.00 & 1.17 & 5.13 & 8.76  & 11.42 &    &   \\
		\midrule
		\multirow{6}{*}{FRI}  & \multirow{3}{*}{R1} & $i$     & 0 & 4    & 3    & 2    & 0     &       &    & \multirow{3}{*}{23.41 min}   \\
		\cmidrule{3-8}
		&                     & $w_{i4}$   &   & 4.47 & 2.34 & 4.86 &       &       &    &   \\
		&                     & $v_{g_{i}i14}$ &   & 0.00 & 4.47 & 6.81 & 11.67 &       &    &   \\
		\cmidrule{2-11}
		& \multirow{3}{*}{R2} & $i$     & 0 & 9    & 5    & 8    & 0     &       &    & \multirow{3}{*}{22.67 min}   \\
		\cmidrule{3-8}
		&                     & $w_{i4}$   &   & 3.16 & 4.77 & 3.69 &       &       &     &  \\
		&                     & $v_{g_{i}i24}$ &   & 0.00 & 3.16 & 7.93 & 11.62 &       &    &   \\
		\midrule
		\multirow{6}{*}{SAT}  & \multirow{3}{*}{R1} & $i$     & 0 & 10   & 3    & 2    & 1     & 12    & 0  & \multirow{3}{*}{24.26 min}   \\
		\cmidrule{3-10}
		&                     & $w_{i5}$   &   & 2.34 & 1.17 & 1.62 & 3.81  & 2.66  &   &    \\
		&                     & $v_{g_{i}i15}$ &   & 0.00 & 2.34 & 3.51 & 5.13  & 8.94  & 11.60 &  \\
		\cmidrule{2-11}
		& \multirow{3}{*}{R2} & $i$     & 0 & 6    & 11   & 4    & 5     & 9     & 0    & \multirow{3}{*}{29.99 min} \\
		\cmidrule{3-10}
		&                     & $w_{i5}$   &   & 2.42 & 4.00 & 1.49 & 1.59  & 1.58  &   &    \\
		&                     & $v_{g_{i}i25}$ &   & 0.00 & 2.42 & 6.42 & 7.91  & 9.50  & 11.08	&  \\
		\botrule
	\end{tabular}
\end{table}

\begin{figure}[]
	\centering
	\includegraphics{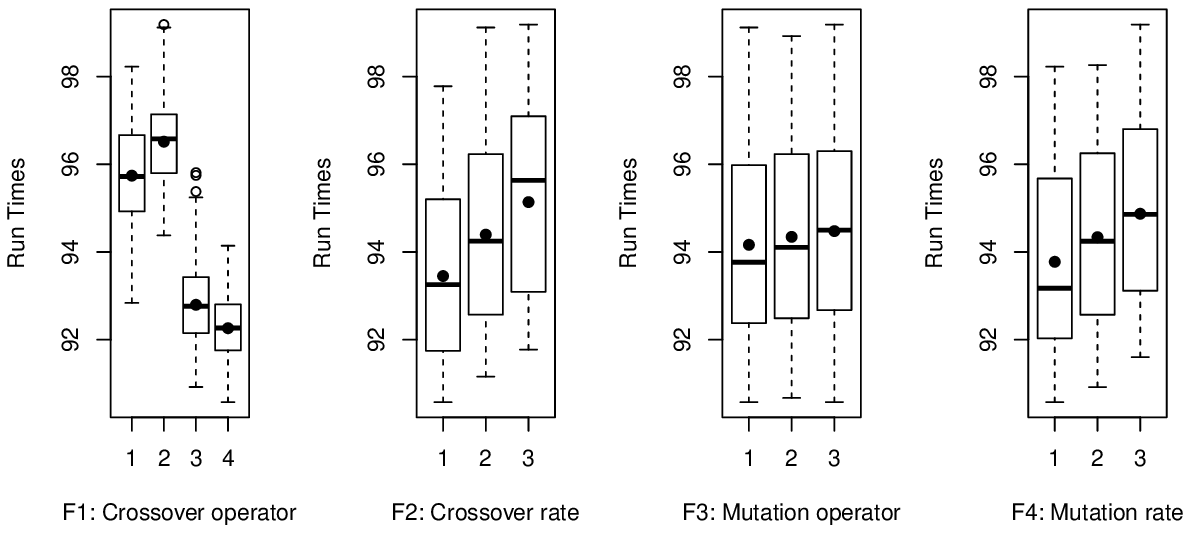}
	\caption{Routes performed by the collection truck each day of the time horizon.}\label{fig:fig5}
\end{figure}

Finally, the routes performed by the collection truck are presented in Figure~\ref{fig:fig5}. Routes are performed from Monday to Saturday (on Sundays no collection is performed). As was represented in Fig.~\ref{fig:fig4}, on Mondays, Thursdays, Fridays and Saturdays, two routes are performed. However, on Wednesday and Thursday only one route is performed. Due to the cycle feature of the proposed PVRP model, this weekly schedule is repeated every week.

\end{appendices}


%\bibliographystyle{}
\bibliography{sn-bibliography}% common bib file


\begin{thebibliography}{39}
\ifx \bisbn   \undefined \def \bisbn  #1{ISBN #1}\fi
\ifx \binits  \undefined \def \binits#1{#1}\fi
\ifx \bauthor  \undefined \def \bauthor#1{#1}\fi
\ifx \batitle  \undefined \def \batitle#1{#1}\fi
\ifx \bjtitle  \undefined \def \bjtitle#1{#1}\fi
\ifx \bvolume  \undefined \def \bvolume#1{\textbf{#1}}\fi
\ifx \byear  \undefined \def \byear#1{#1}\fi
\ifx \bissue  \undefined \def \bissue#1{#1}\fi
\ifx \bfpage  \undefined \def \bfpage#1{#1}\fi
\ifx \blpage  \undefined \def \blpage #1{#1}\fi
\ifx \burl  \undefined \def \burl#1{\textsf{#1}}\fi
\ifx \doiurl  \undefined \def \doiurl#1{\url{https://doi.org/#1}}\fi
\ifx \betal  \undefined \def \betal{\textit{et al.}}\fi
\ifx \binstitute  \undefined \def \binstitute#1{#1}\fi
\ifx \binstitutionaled  \undefined \def \binstitutionaled#1{#1}\fi
\ifx \bctitle  \undefined \def \bctitle#1{#1}\fi
\ifx \beditor  \undefined \def \beditor#1{#1}\fi
\ifx \bpublisher  \undefined \def \bpublisher#1{#1}\fi
\ifx \bbtitle  \undefined \def \bbtitle#1{#1}\fi
\ifx \bedition  \undefined \def \bedition#1{#1}\fi
\ifx \bseriesno  \undefined \def \bseriesno#1{#1}\fi
\ifx \blocation  \undefined \def \blocation#1{#1}\fi
\ifx \bsertitle  \undefined \def \bsertitle#1{#1}\fi
\ifx \bsnm \undefined \def \bsnm#1{#1}\fi
\ifx \bsuffix \undefined \def \bsuffix#1{#1}\fi
\ifx \bparticle \undefined \def \bparticle#1{#1}\fi
\ifx \barticle \undefined \def \barticle#1{#1}\fi
\bibcommenthead
\ifx \bconfdate \undefined \def \bconfdate #1{#1}\fi
\ifx \botherref \undefined \def \botherref #1{#1}\fi
\ifx \url \undefined \def \url#1{\textsf{#1}}\fi
\ifx \bchapter \undefined \def \bchapter#1{#1}\fi
\ifx \bbook \undefined \def \bbook#1{#1}\fi
\ifx \bcomment \undefined \def \bcomment#1{#1}\fi
\ifx \oauthor \undefined \def \oauthor#1{#1}\fi
\ifx \citeauthoryear \undefined \def \citeauthoryear#1{#1}\fi
\ifx \endbibitem  \undefined \def \endbibitem {}\fi
\ifx \bconflocation  \undefined \def \bconflocation#1{#1}\fi
\ifx \arxivurl  \undefined \def \arxivurl#1{\textsf{#1}}\fi
\csname PreBibitemsHook\endcsname

%%% 1
\bibitem[\protect\citeauthoryear{Archetti and Ljubi{\'c}}{2022}]{archetti2022comparison}
\begin{barticle}
\bauthor{\bsnm{Archetti}, \binits{C.}},
\bauthor{\bsnm{Ljubi{\'c}}, \binits{I.}}:
\batitle{Comparison of formulations for the inventory routing problem}.
\bjtitle{European Journal of Operational Research}
\bvolume{303}(\bissue{3}),
\bfpage{997}--\blpage{1008}
(\byear{2022})
\doiurl{10.1016/j.ejor.2021.12.051}
\end{barticle}
\endbibitem

%%% 2
\bibitem[\protect\citeauthoryear{{Asociación Latinoamericana de Logística}}{sept-2023}]{alalog2023}
\begin{botherref}
\oauthor{\bsnm{{Asociación Latinoamericana de Logística}}}:
Indicadores Costos Logísticos países miembros de {ALALOG}.
Available in \url{https://www.alalog.org/es/studies} %(Accessed: sep-2023).
(sept-2023)
\end{botherref}
\endbibitem

%%% 3
\bibitem[\protect\citeauthoryear{Beltrami and Bodin}{1974}]{beltrami74}
\begin{barticle}
\bauthor{\bsnm{Beltrami}, \binits{E.}},
\bauthor{\bsnm{Bodin}, \binits{L.}}:
\batitle{Networks and vehicle routing for municipal waste collection}.
\bjtitle{Networks}
\bvolume{4},
\bfpage{65}--\blpage{94}
(\byear{1974})
\end{barticle}
\endbibitem

%%% 4
\bibitem[\protect\citeauthoryear{Brogaard and Christensen}{2012}]{brogaard2012quantifying}
\begin{barticle}
\bauthor{\bsnm{Brogaard}, \binits{L.}},
\bauthor{\bsnm{Christensen}, \binits{T.}}:
\batitle{Quantifying capital goods for collection and transport of waste}.
\bjtitle{Waste management \& research}
\bvolume{30}(\bissue{12}),
\bfpage{1243}--\blpage{1250}
(\byear{2012})
\doiurl{10.1177/0734242X12462279}
\end{barticle}
\endbibitem

%%% 5
\bibitem[\protect\citeauthoryear{Bynum et~al.}{2021}]{bynum2021pyomo}
\begin{bbook}
\bauthor{\bsnm{Bynum}, \binits{M.}},
\bauthor{\bsnm{Hackebeil}, \binits{G.}},
\bauthor{\bsnm{Hart}, \binits{W.}},
\bauthor{\bsnm{Laird}, \binits{C.}},
\bauthor{\bsnm{Nicholson}, \binits{B.}},
\bauthor{\bsnm{Siirola}, \binits{J.}},
\bauthor{\bsnm{Watson}, \binits{J.}},
\bauthor{\bsnm{Woodruff}, \binits{D.}}:
\bbtitle{Pyomo-optimization Modeling in Python}
vol. \bseriesno{67}.
\bpublisher{Springer}, \blocation{???}
(\byear{2021})
\end{bbook}
\endbibitem

%%% 6
\bibitem[\protect\citeauthoryear{Carlos et~al.}{2019}]{carlos2019influence}
\begin{barticle}
\bauthor{\bsnm{Carlos}, \binits{M.}},
\bauthor{\bsnm{Gallardo}, \binits{A.}},
\bauthor{\bsnm{Edo}, \binits{N.}},
\bauthor{\bsnm{Abaso}, \binits{J.}}:
\batitle{Influence of the municipal solid waste collection system on the time spent at a collection point: {A} case study}.
\bjtitle{Sustainability}
\bvolume{11}(\bissue{22}),
\bfpage{6481}
(\byear{2019})
\doiurl{10.3390/su11226481}
\end{barticle}
\endbibitem

%%% 7
\bibitem[\protect\citeauthoryear{Cavallin et~al.}{2020}]{cavallin2020application}
\begin{barticle}
\bauthor{\bsnm{Cavallin}, \binits{A.}},
\bauthor{\bsnm{Rossit}, \binits{D.}},
\bauthor{\bsnm{Herran}, \binits{V.}},
\bauthor{\bsnm{Rossit}, \binits{D.}},
\bauthor{\bsnm{Frutos}, \binits{M.}}:
\batitle{Application of a methodology to design a municipal waste pre-collection network in real scenarios}.
\bjtitle{Waste Management \& Research}
\bvolume{38}(\bissue{1\_suppl}),
\bfpage{117}--\blpage{129}
(\byear{2020})
\doiurl{10.1177/0734242X19894630}
\end{barticle}
\endbibitem

%%% 8
\bibitem[\protect\citeauthoryear{Cornu{\'e}jols et~al.}{1991}]{cornuejols1991comparison}
\begin{barticle}
\bauthor{\bsnm{Cornu{\'e}jols}, \binits{G.}},
\bauthor{\bsnm{Sridharan}, \binits{R.}},
\bauthor{\bsnm{Thizy}, \binits{J.}}:
\batitle{A comparison of heuristics and relaxations for the capacitated plant location problem}.
\bjtitle{European journal of operational research}
\bvolume{50}(\bissue{3}),
\bfpage{280}--\blpage{297}
(\byear{1991})
\end{barticle}
\endbibitem

%%% 9
\bibitem[\protect\citeauthoryear{Cubillos and Wøhlk}{2020}]{Cubillos2020}
\begin{barticle}
\bauthor{\bsnm{Cubillos}, \binits{M.}},
\bauthor{\bsnm{Wøhlk}, \binits{S.}}:
\batitle{Solution of the maximal covering tour problem for locating recycling drop-off stations}.
\bjtitle{Journal of the Operational Research Society}
\bvolume{72}(\bissue{8}),
\bfpage{1898}--\blpage{1913}
(\byear{2020})
\doiurl{10.1080/01605682.2020.1746701}
\end{barticle}
\endbibitem

%%% 10
\bibitem[\protect\citeauthoryear{Davis}{1985}]{davis85}
\begin{bchapter}
\bauthor{\bsnm{Davis}, \binits{L.}}:
\bctitle{Applying adaptive algorithms to epistatic domains}.
In: \bbtitle{International Joint Conferences on Artificial Intelligence (IJCAI)},
pp. \bfpage{162}--\blpage{164}
(\byear{1985})
\end{bchapter}
\endbibitem

%%% 11
\bibitem[\protect\citeauthoryear{D'Onza et~al.}{2016}]{d2016full}
\begin{barticle}
\bauthor{\bsnm{D'Onza}, \binits{G.}},
\bauthor{\bsnm{Greco}, \binits{G.}},
\bauthor{\bsnm{Allegrini}, \binits{M.}}:
\batitle{Full cost accounting in the analysis of separated waste collection efficiency: {A} methodological proposal}.
\bjtitle{Journal of environmental management}
\bvolume{167},
\bfpage{59}--\blpage{65}
(\byear{2016})
\doiurl{10.1016/j.jenvman.2015.09.002}
\end{barticle}
\endbibitem

%%% 12
\bibitem[\protect\citeauthoryear{Giraudoux et~al.}{2023}]{giraudoux23}
\begin{botherref}
\oauthor{\bsnm{Giraudoux}, \binits{P.}},
\oauthor{\bsnm{Antonietti}, \binits{J.}},
\oauthor{\bsnm{Beale}, \binits{C.}},
\oauthor{\bsnm{Groemping}, \binits{U.}},
\oauthor{\bsnm{Lancelot}, \binits{R.}},
\oauthor{\bsnm{Pleydell}, \binits{D.}},
\oauthor{\bsnm{Treglia}, \binits{M.}}:
{pgirmess}: {S}patial Analysis and Data Mining for Field Ecologists.
(2023).
R package version 2.0.2.
\url{https://www.R-project.org/package=pgirmess}
\end{botherref}
\endbibitem

%%% 13
\bibitem[\protect\citeauthoryear{Goldberg and Lingle}{2014}]{goldberg14}
\begin{bchapter}
\bauthor{\bsnm{Goldberg}, \binits{D.}},
\bauthor{\bsnm{Lingle}, \binits{R.}}:
\bctitle{Alleles, {L}oci and the traveling salesman problem}.
In: \bbtitle{First International Conference on Genetic Algorithms and Their Applications},
pp. \bfpage{154}--\blpage{159}
(\byear{2014}).
\bcomment{Psychology Press}
\end{bchapter}
\endbibitem

%%% 14
\bibitem[\protect\citeauthoryear{Glover}{1975}]{glover1975improved}
\begin{barticle}
\bauthor{\bsnm{Glover}, \binits{F.}}:
\batitle{Improved linear integer programming formulations of nonlinear integer problems}.
\bjtitle{Management science}
\bvolume{22}(\bissue{4}),
\bfpage{455}--\blpage{460}
(\byear{1975})
\end{barticle}
\endbibitem

%%% 15
\bibitem[\protect\citeauthoryear{Gurobi~Optimization}{2023}]{gurobi}
\begin{botherref}
\oauthor{\bsnm{Gurobi~Optimization}, \binits{L.}}:
Gurobi Optimizer Reference Manual
(2023).
\url{https://www.gurobi.com}
\end{botherref}
\endbibitem

%%% 16
\bibitem[\protect\citeauthoryear{Gl{\"a}ser and St{\"u}cken}{2021}]{glaser2021introduction}
\begin{barticle}
\bauthor{\bsnm{Gl{\"a}ser}, \binits{S.}},
\bauthor{\bsnm{St{\"u}cken}, \binits{M.}}:
\batitle{Introduction of an underground waste container system--model and solution approaches}.
\bjtitle{European Journal of Operational Research}
\bvolume{295}(\bissue{2}),
\bfpage{675}--\blpage{689}
(\byear{2021})
\end{barticle}
\endbibitem

%%% 17
\bibitem[\protect\citeauthoryear{Hemmelmayr et~al.}{2014}]{hemmelmayr2014models}
\begin{barticle}
\bauthor{\bsnm{Hemmelmayr}, \binits{V.}},
\bauthor{\bsnm{Doerner}, \binits{K.}},
\bauthor{\bsnm{Hartl}, \binits{R.}},
\bauthor{\bsnm{Vigo}, \binits{D.}}:
\batitle{Models and algorithms for the integrated planning of bin allocation and vehicle routing in solid waste management}.
\bjtitle{Transportation Science}
\bvolume{48}(\bissue{1}),
\bfpage{103}--\blpage{120}
(\byear{2014})
\doiurl{10.1287/trsc.2013.0459}
\end{barticle}
\endbibitem

%%% 18
\bibitem[\protect\citeauthoryear{Hussain et~al.}{2017}]{hussain17}
\begin{botherref}
\oauthor{\bsnm{Hussain}, \binits{A.}},
\oauthor{\bsnm{Muhammad}, \binits{Y.}},
\oauthor{\bsnm{Nauman}, \binits{M.}},
\oauthor{\bsnm{Hussain}, \binits{I.}},
\oauthor{\bsnm{Mohamd}, \binits{A.}},
\oauthor{\bsnm{Gani}, \binits{S.}}:
Genetic algorithm for traveling salesman problem with modified cycle crossover operator.
Computational Intelligence and Neuroscience,
1--7
(2017)
\doiurl{10.1155/2017/7430125}
\end{botherref}
\endbibitem

%%% 19
\bibitem[\protect\citeauthoryear{Holland}{1992}]{holland92}
\begin{bbook}
\bauthor{\bsnm{Holland}, \binits{J.}}:
\bbtitle{Adaptation in Natural and Artificial Systems}.
\bpublisher{MIT Press},
\blocation{Cambridge}
(\byear{1992})
\end{bbook}
\endbibitem

%%% 20
\bibitem[\protect\citeauthoryear{Hemmelmayr et~al.}{2017}]{hemmelmayr2017periodic}
\begin{barticle}
\bauthor{\bsnm{Hemmelmayr}, \binits{V.}},
\bauthor{\bsnm{Smilowitz}, \binits{K.}},
\bauthor{\bsnm{Torre}, \binits{L.}}:
\batitle{A periodic location routing problem for collaborative recycling}.
\bjtitle{IISE Transactions}
\bvolume{49}(\bissue{4}),
\bfpage{414}--\blpage{428}
(\byear{2017})
\doiurl{10.1080/24725854.2016.1267882}
\end{barticle}
\endbibitem

%%% 21
\bibitem[\protect\citeauthoryear{Hollander and Wolfe}{1973}]{hollander73}
\begin{bbook}
\bauthor{\bsnm{Hollander}, \binits{M.}},
\bauthor{\bsnm{Wolfe}, \binits{D.}}:
\bbtitle{Nonparametric Statistical Methods}.
\bpublisher{John Wiley \& Sons},
\blocation{New Jersey}
(\byear{1973})
\end{bbook}
\endbibitem

%%% 22
\bibitem[\protect\citeauthoryear{Han et~al.}{2024}]{Han2024Optimizing}
\begin{barticle}
\bauthor{\bsnm{Han}, \binits{J.}},
\bauthor{\bsnm{Zhang}, \binits{J.}},
\bauthor{\bsnm{Guo}, \binits{H.}},
\bauthor{\bsnm{Zhang}, \binits{N.}}:
\batitle{Optimizing location-routing and demand allocation in the household waste collection system using a branch-and-price algorithm}.
\bjtitle{European Journal of Operational Research}
\bvolume{316}(\bissue{3}),
\bfpage{958}--\blpage{975}
(\byear{2024})
\end{barticle}
\endbibitem

%%% 23
\bibitem[\protect\citeauthoryear{Maalouf and Agamuthu}{2023}]{maalouf2023waste}
\begin{barticle}
\bauthor{\bsnm{Maalouf}, \binits{A.}},
\bauthor{\bsnm{Agamuthu}, \binits{P.}}:
\batitle{Waste management evolution in the last five decades in developing countries--{A} review}.
\bjtitle{Waste Management \& Research}
\bvolume{41}(\bissue{9}),
\bfpage{1420}--\blpage{1434}
(\byear{2023})
\doiurl{10.1177/0734242X231160099}
\end{barticle}
\endbibitem

%%% 24
\bibitem[\protect\citeauthoryear{Mah{\'e}o et~al.}{2020}]{maheo2020benders}
\begin{botherref}
\oauthor{\bsnm{Mah{\'e}o}, \binits{A.}},
\oauthor{\bsnm{Rossit}, \binits{D.}},
\oauthor{\bsnm{Kilby}, \binits{P.}}:
A {B}enders decomposition approach for an integrated bin allocation and vehicle routing problem in municipal waste management.
Communications in Computer and Information Science,
3--18
(2020)
\doiurl{10.1007/978-3-030-76310-7_1}
\end{botherref}
\endbibitem

%%% 25
\bibitem[\protect\citeauthoryear{Mah{\'e}o et~al.}{2023}]{maheo2023solving}
\begin{barticle}
\bauthor{\bsnm{Mah{\'e}o}, \binits{A.}},
\bauthor{\bsnm{Rossit}, \binits{D.}},
\bauthor{\bsnm{Kilby}, \binits{P.}}:
\batitle{Solving the integrated bin allocation and collection routing problem for municipal solid waste: {A} {B}enders decomposition approach}.
\bjtitle{Annals of Operations Research}
\bvolume{322}(\bissue{1}),
\bfpage{441}--\blpage{465}
(\byear{2023})
\doiurl{10.1007/s10479-022-04918-7}
\end{barticle}
\endbibitem

%%% 26
\bibitem[\protect\citeauthoryear{Niu et~al.}{2024}]{niu2024multi}
\begin{barticle}
\bauthor{\bsnm{Niu}, \binits{Y.}},
\bauthor{\bsnm{Xu}, \binits{C.}},
\bauthor{\bsnm{Liao}, \binits{S.}},
\bauthor{\bsnm{Zhang}, \binits{S.}},
\bauthor{\bsnm{Xiao}, \binits{J.}}:
\batitle{Multi-objective location-routing optimization based on machine learning for green municipal waste management}.
\bjtitle{Waste Management}
\bvolume{181},
\bfpage{157}--\blpage{167}
(\byear{2024})
\doiurl{10.1016/j.wasman.2024.04.001}
\end{barticle}
\endbibitem

%%% 27
\bibitem[\protect\citeauthoryear{Owusu et~al.}{2019}]{Owusu-Nimo2019}
\begin{botherref}
\oauthor{\bsnm{Owusu}, \binits{F.}},
\oauthor{\bsnm{Oduro}, \binits{S.}},
\oauthor{\bsnm{Essandoh}, \binits{H.}},
\oauthor{\bsnm{Wayo}, \binits{F.}},
\oauthor{\bsnm{Shamudeen}, \binits{M.}}:
Characteristics and management of landfill solid waste in {Kumasi, Ghana}.
Scientific African
\textbf{3}(e00052)
(2019)
\doiurl{10.1016/j.sciaf.2019.e0 0 052}
\end{botherref}
\endbibitem

%%% 28
\bibitem[\protect\citeauthoryear{Oliver et~al.}{1987}]{oliver87}
\begin{bchapter}
\bauthor{\bsnm{Oliver}, \binits{I.}},
\bauthor{\bsnm{Smith}, \binits{D.}},
\bauthor{\bsnm{Holland}, \binits{J.}}:
\bctitle{A study of permutation crossover operators on the traveling salesman problem}.
In: \bbtitle{Second International Conference on Genetic Algorithms on Genetic Algorithms and Their Application},
pp. \bfpage{224}--\blpage{230}
(\byear{1987}).
\bcomment{Psychology Press}
\end{bchapter}
\endbibitem

%%% 29
\bibitem[\protect\citeauthoryear{{Posit team}}{2023}]{RStudio2023}
\begin{bbook}
\bauthor{\bsnm{{Posit team}}}:
\bbtitle{RStudio: {I}ntegrated Development Environment For {R}}.
\bpublisher{Posit Software, PBC},
\blocation{Boston, MA}
(\byear{2023}).
\bcomment{Posit Software, PBC}.
\burl{http://www.posit.co/}
\end{bbook}
\endbibitem

%%% 30
\bibitem[\protect\citeauthoryear{{R Core Team}}{2022}]{RCoreTeam2022}
\begin{bbook}
\bauthor{\bsnm{{R Core Team}}}:
\bbtitle{R: {A} Language and Environment for Statistical Computing}.
\bpublisher{R Foundation for Statistical Computing},
\blocation{Vienna, Austria}
(\byear{2022}).
\bcomment{R Foundation for Statistical Computing}.
\burl{https://www.R-project.org/}
\end{bbook}
\endbibitem

%%% 31
\bibitem[\protect\citeauthoryear{Rossit et~al.}{2024}]{Rossit2024}
\begin{bchapter}
\bauthor{\bsnm{Rossit}, \binits{D.}},
\bauthor{\bsnm{Gonz{\'a}lez~Land{\'\i}n}, \binits{B.}},
\bauthor{\bsnm{Frutos}, \binits{M.}},
\bauthor{\bsnm{M{\'e}ndez~Babey}, \binits{M.}}:
\bctitle{An allocation-routing problem in waste management planning: {E}xact and heuristic resolution approaches}.
In: \beditor{\bsnm{Nesmachnow}, \binits{S.}},
\beditor{\bsnm{Hernández~Callejo}, \binits{L.}} (eds.)
\bbtitle{Smart Cities. ICSC-CITIES 2023. Communications in Computer and Information Science},
vol. \bseriesno{1938},
pp. \bfpage{92}--\blpage{107}.
\bpublisher{Springer},
\blocation{Cham.}
(\byear{2024}).
\doiurl{10.1007/978-3-031-52517-9_7}
\end{bchapter}
\endbibitem

%%% 32
\bibitem[\protect\citeauthoryear{Roy et~al.}{2022}]{roy2022iot}
\begin{barticle}
\bauthor{\bsnm{Roy}, \binits{A.}},
\bauthor{\bsnm{Manna}, \binits{A.}},
\bauthor{\bsnm{Kim}, \binits{J.}},
\bauthor{\bsnm{Moon}, \binits{I.}}:
\batitle{Io{T}-based smart bin allocation and vehicle routing in solid waste management: {A case study in South Korea}}.
\bjtitle{Computers \& Industrial Engineering}
\bvolume{171},
\bfpage{108457}
(\byear{2022})
\doiurl{10.1016/j.cie.2022.108457}
\end{barticle}
\endbibitem

%%% 33
\bibitem[\protect\citeauthoryear{Rossit and Nesmachnow}{2022}]{rossit2022waste}
\begin{barticle}
\bauthor{\bsnm{Rossit}, \binits{D.}},
\bauthor{\bsnm{Nesmachnow}, \binits{S.}}:
\batitle{Waste bins location problem: {A} review of recent advances in the storage stage of the municipal solid waste reverse logistic chain}.
\bjtitle{Journal of Cleaner Production}
\bvolume{342},
\bfpage{130793}
(\byear{2022})
\doiurl{10.1016/j.jclepro.2022.130793}
\end{barticle}
\endbibitem

%%% 34
\bibitem[\protect\citeauthoryear{Rossit et~al.}{2023}]{rossit2023systemwaste}
\begin{bchapter}
\bauthor{\bsnm{Rossit}, \binits{D.}},
\bauthor{\bsnm{Nesmachnow}, \binits{S.}},
\bauthor{\bsnm{Cavallin}, \binits{A.}}:
\bctitle{Municipal solid waste management systems: {A}pplication of {SWOT} methodology to analyze an {A}rgentinean case study}.
In: \bbtitle{VI Ibero-American Congress of Smart Cities ICSC-CITIES 2023},
\bconflocation{Mexico City},
pp. \bfpage{668}--\blpage{680}
(\byear{2023})
\end{bchapter}
\endbibitem

%%% 35
\bibitem[\protect\citeauthoryear{Rossi et~al.}{2022}]{rossi2022comparison}
\begin{barticle}
\bauthor{\bsnm{Rossi}, \binits{M.}},
\bauthor{\bsnm{Papetti}, \binits{A.}},
\bauthor{\bsnm{Germani}, \binits{M.}}:
\batitle{A comparison of different waste collection methods: {E}nvironmental impacts and occupational risks}.
\bjtitle{Journal of Cleaner Production}
\bvolume{368},
\bfpage{133145}
(\byear{2022})
\doiurl{10.1016/j.jclepro.2022.133145}
\end{barticle}
\endbibitem

%%% 36
\bibitem[\protect\citeauthoryear{Rossit et~al.}{2017}]{rossit2017application}
\begin{barticle}
\bauthor{\bsnm{Rossit}, \binits{D.}},
\bauthor{\bsnm{Tohm{\'e}}, \binits{F.}},
\bauthor{\bsnm{Frutos}, \binits{M.}},
\bauthor{\bsnm{Broz}, \binits{D.}}:
\batitle{An application of the augmented $\varepsilon$-constraint method to design a municipal sorted waste collection system}.
\bjtitle{Decision Science Letters}
\bvolume{6}(\bissue{4}),
\bfpage{323}--\blpage{336}
(\byear{2017})
\doiurl{10.5267/j.dsl.2017.3.001}
\end{barticle}
\endbibitem

%%% 37
\bibitem[\protect\citeauthoryear{Rodriguez and Vecchietti}{2013}]{rodriguez2013comparative}
\begin{barticle}
\bauthor{\bsnm{Rodriguez}, \binits{M.}},
\bauthor{\bsnm{Vecchietti}, \binits{A.}}:
\batitle{A comparative assessment of linearization methods for bilinear models}.
\bjtitle{Computers \& Chemical Engineering}
\bvolume{48},
\bfpage{218}--\blpage{233}
(\byear{2013})
\doiurl{10.1016/j.compchemeng.2012.09.011}
\end{barticle}
\endbibitem

%%% 38
\bibitem[\protect\citeauthoryear{Siegel and Castellan}{1988}]{siegel88}
\begin{bbook}
\bauthor{\bsnm{Siegel}, \binits{S.}},
\bauthor{\bsnm{Castellan}, \binits{N.}}:
\bbtitle{Non Parametric Statistics for the Behavioural Sciences}.
\bpublisher{MacGraw Hill},
\blocation{New York}
(\byear{1988})
\end{bbook}
\endbibitem

%%% 39
\bibitem[\protect\citeauthoryear{Toth and Vigo}{2014}]{toth2014vehicle}
\begin{bbook}
\bauthor{\bsnm{Toth}, \binits{P.}},
\bauthor{\bsnm{Vigo}, \binits{D.}}:
\bbtitle{Vehicle Routing: {P}roblems, Methods, and Applications}.
\bpublisher{Society for Industrial and Applied Mathematics},
\blocation{Philadelphia}
(\byear{2014}).
\doiurl{10.1137/1.9781611973594}
\end{bbook}
\endbibitem

\end{thebibliography}
%% \input sn-article.bbl

\end{document}